\documentclass[11pt]{article}

\usepackage{epsfig}
\usepackage[latin1]{inputenc}
\usepackage{amsmath}
\usepackage{amssymb}
\usepackage{psfrag} 
\usepackage{version} 
\usepackage{dsfont} 

\def \bye {\end{document}}

\let \article =\relax 

\def \SECTION #1#2{\section {#2\label {s #1}}} 
\def \TSECTION #1{\section* {#1}}
\def \SUBSECTION #1#2{\subsection {#2\label {ss #1}}} 
\def \SSUBSECTION #1#2{\subsubsection {#2\label {sss #1}}} 
\def \sectionref #1{\ref {s #1}} 
\def \subsectionref #1{\ref {ss #1}}
\def \ssubsectionref #1{\ref {sss #1}}

\newtheorem{theorem}{Theorem}[section]
\def \THEOREM #1{\begin {theorem}\label {t #1}} 
\def \ENDTHEOREM {\end {theorem}} 
\def \PROPOSITION #1{\begin {proposition}\label {t #1}} 
\def \ENDPROPOSITION {\end {proposition}} 
\newtheorem{lemma}[theorem]{Lemma} 
\def \LEMMA #1{\begin {lemma}\label {t #1}} 
\def \ENDLEMMA {\end {lemma}} 
\newtheorem{corollary}[theorem]{Corollary} 
\def \COROLLARY #1{\begin {corollary}\label {t #1}} 
\def \ENDCOROLLARY {\end {corollary}} 
\def \statementref #1{\ref {t #1}} 
\newtheorem{definition}[theorem]{Definition} 
\def \DEFINITION #1{\begin {definition}\label {t #1}} 
\def \ENDDEFINITION {\end {definition}} 
\def \statementref #1{\ref {t #1}} 

\newenvironment {remark}[1]{\bigbreak \noindent {\sl #1\space }\ignorespaces }{\bigbreak } 
\def \REMARK{\begin{remark}{Remark.}} 
\def \ENDREMARK {\end {remark}} 
\newenvironment {proof}[1]{\bigbreak \noindent {\bf #1\space }\ignorespaces }{\bigbreak } 
\def \PROOF{\begin{proof}{Proof.}} 
\def \ENDPROOF {\end {proof}} 

\def \eqref #1{(\ref {#1})}

\def \square {\hbox {$\sqcup $\llap {$\sqcap $}}} 
\def \endsquare {\ifmmode \kern 0.5em\square \else \discretionary {} {\hbox to\hsize {\hfill \square }\null }{\kern 0.5em\square }\fi } 

\def \seventeenbf {\huge \bf } 

\def \newline {\ifhmode \hfill \break \else {\everypar ={}\noindent }\fi } 

\long \def \centerlines #1{{\center #1\par }} 

\def \BIB {\begin {thebibliography}{2}} 
\def \ENDBIB {\end {thebibliography}} 
\let \bib =\bibitem 
\let \bibref =\cite 

\def \cD {{\cal D}} 
\let \al =\alpha 
\let \be =\beta

\let \ga =\gamma 
\let \la =\lambda 
\let \Om =\Omega 
\let \na = \nabla 
\let \pa = \partial 
\let \vep =\varepsilon 
\let \vrh =\rho 
\def\vph{\varphi}

\def\dx{\partial_x}

\def\dz{\partial_z}

\def\dZ{\partial_Z}

\def\lap{\Delta}

\def\dsp{\displaystyle}

\def\R{\mathbb{R}}

\def\T{\mathbb{T}}

\def\f{{\bf f}}

\def\m{{\bf m}}
\def\n{{\bf n}}
\def\0{{\bf 0}}

\def\u{{\bf u}}
\def\tu{{\widetilde{\bf u}}}

\def\U{{\bf U}}

\def\x{{\bf x}}

\def\eps{\varepsilon}

\def\fy{\varphi}

\def\bsigma{\boldsymbol{\sigma}}

\def\div{\mbox{\,{\textrm{div}}}}
\newcommand{\nfrac}[2]{#1/#2}

\begin{document}


\article
\parindent=0pt
\def\div{{\rm div}}
\def\rot{{\rm curl}}


\centerlines{\seventeenbf Existence result for stationary compressible fluids and asymptotic behavior in thin films}

\centerlines{
Laurent {\sc Chupin}
  \footnote{Université de Lyon - INSA de Lyon - Pôle de Mathématiques - CNRS, UMR5208, Institut Camille Jordan - 21~av. Jean Capelle, 69621 Villeurbanne Cedex, France.
\newline  e-mail: laurent.chupin@insa-lyon.fr},
R\'emy {\sc Sart}
  \footnote{Institut de Math\'ematiques et Mod\'elisation de Montpellier - Place Eug\`ene Bataillon
34095 Montpellier
\newline  e-mail: rsart@math.univ-montp2.fr} 
}

\bigskip
\bigskip


\begin{abstract}
In this paper, we are first interested in the compressible Navier-Stokes equations with density-dependent viscosities in bounded domains with non-homogeneous Dirichlet conditions.
We study the wellposedness of such models with non-constant coefficients in non-stationary and stationary cases.
We apply the last result in thin domains context, justifying the compressible Reynolds equations.
\end{abstract}

\smallskip
\smallskip
\noindent {\bf Keywords}:
{\small Navier-Stokes equations, compressible fluids, density dependent viscosities, stationary, steady-state, thin films, lubrication, Reynolds equation}

\noindent {\bf AMS subject classification}: 35Q30, 76A20, 76D05, 76D08, 76N10


\TSECTION{Introduction}

The Reynolds equation is a linear equation describing the evolution of the pressure in mechanisms of lubrication.
More precisly, it is used to calculate the pressure distribution in a thin layer of lubricant between two surfaces.
It was proposed by O.~Reynolds in 1886, see~\cite{Reynolds}.
This equation is very much used in mechanics, for instance to describe the process of lubrication of magnetic hard discs.
In the same way, air flow between the two surfaces constituting a rigid disk assembly (flying head and magnetic storage surface) is frequently modeled using the Reynolds equation.
\par\noindent
It was proven, first in 1986 by G.~Bayada and M.~Chambat, see~\cite{Bayada-Chambat}, that the Reynolds equation is an approximation of the Stokes equations in thin cases.
This proof was formulated by taking as initial equations the incompressible Stokes model.
Since, many other works (see~\cite{Nazarov} and the cited references) has made it possible to refine the first result by giving errors of the approximations between the Stokes model and the Reynolds incompressible model.
\par\noindent
In the previous examples of applications the fluids (like air) are clearly compressible fluids.
Within the framework of the compressible fluids, there exist a so called compressible Reynolds equation which is, at least formally, the asymptotic of the Navier-Stokes compressible equations in a thin domain.
Contrary to the incompressible classical Reynolds equation, the compressible Reynolds equation is highly nonlinear and has been a subject of many mechanical studies~\cite{Bhushan,Crone} or of numerical studies~\cite{Arregui,Buscaglia}.
\par\noindent
However the literature about the rigorous justification of these equations in the compressible case is not very important.
It would seem that there is only one result, due to E.~Marusic-Paloka and M.~Starcevic (see~\cite{MarPal-Sta} or more recently~\cite{MarPal-Sta09}), restricted to the case of ideal gases.
The primary reason of this lack of literature certainly comes from the fact that the study of the compressible Navier-Stokes equations is rather difficult.
Recent works of D.~Bresch and B.~Desjardins on these compressible Navier-Stokes equations, see~\cite{BD} for instance, showed that there exists a particular structure to these equations.
\par\noindent
In this article, we adapt these new results to use them and rigorously justify the compressible Reynolds equation for rather general state laws.
\par\noindent
More precisely, we prove the two following results (precise statements are respectively given on page~\pageref{th-stationary} and page~\pageref{th-convergence}):
\begin{theorem}
There exists a steady-state solution to the compressible Navier-Stokes equations in a bounded domain, with Dirichlet boundary conditions.
\end{theorem}

\begin{theorem}
The compressible Reynolds equation is an approximation of the stationary compressible Navier-Stokes equations.
\end{theorem}

The present paper is composed of the following parts:
\begin{itemize}
\item In the first section, we present the notations and the classical compressible Navier-Stokes equations.
We give the assumptions as well as the theorems related to the compressible Navier-Stokes equations (in the non-stationary case and in the stationary case).
\item Section~\ref{sec-2} is devoted to the proof of the existence result for the compressible Navier-Stokes equations in the non stationary case.
\item Section~\ref{sec-3} is devoted to the proof of the existence result for the compressible Navier-Stokes equations in the stationary case.
\item In Section~\ref{sec:thin}, we introduce the lubrication problem in term of thin film flow.
We also anounce the convergence result for the stationary compressible Navier-Stokes equations to the compressible Reynolds equation.
\item Section~\ref{sec-5} is devoted to the proof of this convergence result.
\end{itemize}

\SECTION{1}{Wellposedness of compressible Navier-Stokes equations with\\ density-dependent viscosities}

\SUBSECTION{1.1}{Statement of the problem}

{\bf Compressible Navier-Stokes equations}:
The compressible Navier-Stokes equations describe the evolution of a compressible fluid in a physical domain~$\Omega\subset \R^d$, $d\in \{2,3\}$, {\it via} the conservation equations of the mass and the momentum. They thus couple the velocity~$\u$ of the fluid and its density~$\rho$:
\begin{equation*}
\left\{
\begin{aligned}
& \partial_t \rho + \div (\rho \u) = 0,\\
& \partial_t (\rho \u) + \div (\rho \u \otimes \u) = \div (\bsigma) + \f.
\end{aligned}
\right.
\end{equation*}
To close this system, we must give the force term~$\f$ and the stress tensor~$\bsigma$.
\begin{itemize}
\item[*] The force term~$\f$ allows to represent friction forces or corresponds to a turbulent
drag force. They read
\begin{equation*}
\f = - r_0 \rho |\u| \u,
\end{equation*}
 where $r_0$ is a non negative real coefficient.
\item[*] Finally, we give the rheological law for the stress tensor~$\bsigma$ : the fluid is assumed to be Newtonian, so that there exists two viscosity coefficients (called Lam\'e coefficients) $\mu=\mu(\vrh)$ and $\lambda=\lambda(\vrh)$ such that
\begin{equation*}
  \bsigma = 2\mu D(\u) + (\lambda \div (\u) - p){\rm Id}.
\end{equation*}
In this equation, $D(\u)$ corresponds to the strain tensor (the symmetric part of the velocity gradient tensor), and the pressure~$p$ is determined using a thermodynamic closure law, i.e. an explicit relation $p=p(\vrh)$.
\end{itemize}
The equations in which we will be interested are thus:
\begin{equation}\label{1}
  \partial_t \rho + \div (\rho \u) = 0,
\end{equation}
\begin{equation}\label{2}
  \partial_t (\rho \u) + \div (\rho \u \otimes \u) + \nabla p(\rho) = \div (2\mu(\rho)D(\u)) + \nabla (\lambda(\rho)\div (\u)) - r_0 \rho |\u| \u,
\end{equation}
and also the stationary correponding ones:
\begin{equation}\label{1stat}
  \div (\rho \u) = 0,
\end{equation}
\begin{equation}\label{2stat}
  \div (\rho \u \otimes \u) + \nabla p(\rho) = \div (2\mu(\rho)D(\u)) + \nabla (\lambda(\rho)\div (\u)) - r_0 \rho |\u| \u.
\end{equation}
{\bf Boundary conditions}: The physical boundary conditions which interest us here (see part~\ref{sec:thin}) are of the nonhomogeneous Dirichlet type on the velocity field. For technical reasons, we will also need to impose a condition on the density on the boundary. The conditions will thus be the following ones:
\begin{equation}\label{boundary}
\begin{aligned}
&  \u=\u_b \quad \text{a given function on $\partial \Omega$ such that $\u_b\cdot \n = 0$,}\\
&  \rho=\rho_b \quad \text{constant on each connected component of $\partial \Omega$.}
\end{aligned}
\end{equation}
The biggest part of the published works concerns the whole space case $\Omega=\R^3$, or the periodic case $\Omega=\T^3$ (see for instance~\cite{BD,BDL}).
More recently in~\cite{BDG}, the authors deals with the Dirichlet homogeneous condition or Navier's condition on the velocity field.
The building that we present here draws hard inspiration from this last paper. Particularly, the boundary condition on~$\rho$, being already present in~\cite{BDG}, it is not amazing to find it in a more general case.\\[0.2cm]
{\bf Initial conditions}: In the non-stationary case, it is necessary to give the initial conditions corresponding to the situation at time $t=0$. The physical quantities for which we give information are the density and the momentum:
\begin{equation}\label{initial}
  \rho|_{t=0}=\rho_0
  \quad \text{and}  \quad
  \rho \u|_{t=0}=\m_0.
\end{equation}
We note that, due to the non-penetration condition $\u_b\cdot \n = 0$, integrating with respect to the spatial variable the mass conservation equation~\eqref{1} we obtain
\begin{equation*}
\frac{d}{dt} \Big( \int_\Omega \rho\Big) = 0.
\end{equation*}
Consequently, the quantity $\int_\Omega \rho$ does not depend on time and will be denoted~$M_0$.\\

This system~\eqref{1}-\eqref{2} has been widely studied, starting from the case of constant coefficients $\lambda$, $\mu$ and pressure laws of type $p(\rho)=a\rho^\gamma$ (see notably~\cite{FNP,F,L,L2,NS}). More recently, many studies have focused on density dependent
viscosity coefficients $\lambda=\lambda(\rho)$, $\mu=\mu(\rho)$ in space dimensions~$2$ or~$3$. These studies were originally developed on Korteweg and shallow water models, corresponding to $\gamma = 2$, $\lambda(\rho)=0$ and $\mu(\rho)=\rho$, see~\cite{BD1,BD2,BD,BDG,BDL}. They all rely on a new mathematical entropy (the BD entropy), that has been discovered in its general form in \bibref{BD}. It requires that the
following algebraic relation holds:
\begin{equation*}
  \forall s>0,\ \la(s)=2(s\mu'(s)-\mu(s)).
\end{equation*}
We introduce in the next part, the hypotheses which we shall use later. Obviously, these hypotheses take back principally those of papers named here.

\SUBSECTION{1.2}{Assumptions}


Concerning the viscosity coefficients $\lambda$ and $\mu$, we assume that $\la$ and $\mu$ are respectively~$C^0(\R_+)$ and~$C^1(\R_+)$ and satisfy
\begin{equation}\label{5}
  \forall s>0,\ \la(s)=2(s\mu'(s)-\mu(s)).
\end{equation}
We also suppose that $\mu(0)=0$, that there exists positive constants $c_0,c_1,c_0',c_1',A$, $m>1$ and $\frac{2}{3}<n<1$ such that
\begin{equation}\label{6}
  \forall s\in ]0,A[,\quad c_0\, s^n\leq\mu(s)\leq \frac{1}{c_0} s^n,\quad  c_0'\, s^{n-1}\leq\mu'(s)\leq \frac{1}{c_0'} s^{n-1},
\end{equation}
\begin{equation}\label{7}
  \forall s\in ]A,+\infty[,\quad c_1 s^m\leq\mu(s)\leq \frac{1}{c_1} s^m,\quad c_1' s^{m-1}\leq\mu'(s)\leq \frac{1}{c_1'} s^{m-1}.
\end{equation}
\noindent Finally, we are interested in a pressure term of the following form
\begin{equation}\label{8}
  p(\vrh) = p_h(\vrh) + p_c(\vrh),
\end{equation}
where $p_h(\vrh)=a\vrh^\gamma$ ($a>0$ and $\gamma \geq 1$) corresponds to the classical equation of state whereas $p_c(\vrh)$ is a "cold" component.
We assume that there exists positive constants $c_2$, $c_3$, $\vrh_*$, $\beta$ and~$\alpha\geq 1$ such that
\begin{equation}\label{9}
  \forall \vrh\in ]0,\vrh_*[, \quad \frac{1}{c_2\, \vrh^{\alpha+1}} \leq 
  p'_c(\vrh) \leq \frac{c_2}{\vrh^{\alpha+1}},
\end{equation}
\begin{equation}\label{10}
  \forall \vrh\in ]\vrh_*,+\infty[, \quad -\frac{a \gamma 
  \vrh^{\gamma-1}}{2} \leq p'_c(\vrh) \leq c_3\, \vrh^{\beta-1}.
\end{equation}
Recall that such assumptions were initially introduced in \bibref{BD} in the framework of barotropic flows.
\par
Moreover, in the three-dimensional case, we impose (in fact in the stationary case result)
\begin{equation}\label{condition}
m < \ga+n-\frac{1}{3},
\quad
\beta \leq 2(\ga +n-1),
\end{equation}
and (to control the ``thin domain''-dependency)
\begin{equation}\label{condition+}
m < \al-n+\frac{7}{3}.
\end{equation}

\SUBSECTION{1.3}{Existence results}

\DEFINITION{1}
We shall say that $(\rho, \u)$ is a weak solution of \eqref{1}--\eqref{2} with boundary conditions \eqref{boundary} if it satisfies following regularity properties
\begin{equation*}
\begin{aligned}
& \rho\in L^\infty(0,T;L^\gamma(\Omega)),\quad
\sqrt{\rho}\u \in L^\infty(0,T;L^2(\Omega)),\quad \sqrt{\rho} \nabla \fy(\rho) \in L^\infty(0,T;L^2(\Omega)),\\
& \sqrt{\mu(\rho)}\nabla \u \in L^2((0,T)\times \Omega),\quad
\rho\u^{3} \in L^1((0,T)\times \Omega),\\
& \na\rho^{\frac{\ga+n-1}{2}}\in L^2(0,T;H^1(\Omega)),\quad \na\xi(\rho)^{\frac{n-\al-1}{2}}\in L^2(0,T;H^1(\Omega)),
\end{aligned}
\end{equation*}
where $\xi$ being taken such that $\xi(\rho)=\rho$ for $\rho\leq \rho_*/2$ and $\xi(\rho)=0$ for $\rho\geq \rho_*$, as well as boundary Dirichlet conditions on $\u$ in $L^2(0,T;L^1(\partial \Omega))$, boundary conditions on $\rho$ in $L^2(0,T;L^\infty(\partial \Omega))$, and equations \eqref{1}--\eqref{2} in $\mathcal D'((0,T)\times \Omega)$ for all $T >0$.
\ENDDEFINITION
As usual, we deduce from these regularities and the Navier-Stokes system itself that~$\rho$
and~$\u$ are continuous in time with values in $W^{-1,1}(\Omega)$, which allows to define their initial values. 
%
\begin{theorem}[Non-stationary case]\label{th-non-stationary}
Assume that conditions \eqref{5}--\eqref{10} are satisfied and consider some functions $\vrh_0$ and $\m_0$ such that 
\[
\frac{\m_0^2}{\vrh_0}\in L^1(\Om),\ \ \frac{|\na\mu(\vrh_0)|^2}{\vrh_0}\in L^1(\Om),\ \ Q(\vrh_0)\in L^1(\Om),
\]
where $xQ''(x):=p(x)$.

Then, for all $r_0\in\R^+$, there exists a weak solution of the system 
\begin{eqnarray}
  \partial_t \rho + \div (\rho \u) & = & 0,\nonumber \\
  \partial_t (\rho \u) + \div (\rho \u \otimes \u) + \nabla p(\rho) & = & \div (2\mu(\rho)D(\u)) + \nabla (\lambda(\rho)\div (\u)) - r_0 \rho |\u| \u,\nonumber
\end{eqnarray}
associated with the initial conditions \eqref{initial}, in the sense of the Definition \statementref{1}.
\end{theorem}
\DEFINITION{1stat}
We shall say that $(\rho, \u)$ is a weak solution of \eqref{1stat}--\eqref{2stat} with boundary conditions \eqref{boundary} if it satisfies following regularity properties
\begin{equation*}
\begin{aligned}
& \sqrt{\mu(\rho)}\nabla \u \in L^2(\Omega),\quad
\rho\u^{3} \in L^1(\Omega),\\
& \na\rho^{\frac{\ga+n-1}{2}}\in H^1(\Omega), \quad \na\xi(\rho)^{\frac{n-\al-1}{2}}\in H^1(\Omega),
\end{aligned}
\end{equation*}
where $\xi$ being taken such that $\xi(\rho)=\rho$ for $\rho\leq \rho_*/2$ and $\xi(\rho)=0$ for $\rho\geq \rho_*$,
as well as boundary Dirichlet conditions on~$\u$ in $L^1(\partial \Omega)$, boundary conditions on~$\rho$ in $L^\infty(\partial \Omega)$, and equations \eqref{1stat}--\eqref{2stat} in~$\mathcal D'(\Omega)$.
\ENDDEFINITION

\begin{theorem}[Stationary case]\label{th-stationary}
Assume that conditions \eqref{5}--\eqref{condition} are satisfied.
\par\noindent
Then, for all $r_0\in\R^+_*$, there exists a weak solution of the stationary equations 
\begin{eqnarray}
  \div (\rho \u) & = & 0,\nonumber \\
  \div (\rho \u \otimes \u) + \nabla p(\rho) & = & \div (2\mu(\rho)D(\u)) + \nabla (\lambda(\rho)\div (\u)) - r_0 \rho |\u| \u,\nonumber
\end{eqnarray}
in the sense of the Definition~\statementref{1stat}.
\end{theorem}

\SECTION{2}{Proof of Theorem \ref{th-non-stationary} (non-stationary case)}\label{sec-2}

The main idea is to obtain good energy estimates, notably using the BD entropy, structure discovered by D.~Bresh and B.~Desjardins in~\cite{BD1}.
This will provide enough compactness on a sequence of approximate solutions to pass to the limit and obtain a global weak solution.\par
More precisely, the first step is to obtain suitable a priori bounds on $(\rho,\u)$, and next to consider sequences $(\rho_k,\u_k)$ of uniformly bounded weak solutions constructed from an adapted approximation process. Such sequences may be built by using the regularization scheme given in Section~\ref{sub:regul} (see also~\cite{BD3}). It leads to regular approximate solutions, still preserving physical bounds and the mathematical entropy, uniformly with respect to smoothing parameters.\par
The scheme of the proof will be therefore the following. In Section~\ref{sec:BDstrategy}, we recall the main idea of the Bresch-Desjardins strategy, we then deduce energy estimate and so called BD estimate (Subsection~\ref{sec:estimates}).
In Subsection~\ref{sub:regul}, we give the construction of approximate solutions.

\SUBSECTION{2.1}{Bresch-Desjardins strategy}\label{sec:BDstrategy}

The BD entropy is the dedicated idea to get many wellposedness of non-stationary models with non-constant coefficients, for instance the compressible Navier-Stokes equations (see~\bibref{BD}) and it has been recently enlarged to some neighbour contexts like Shallow-Water (see~\bibref{BD3}) or Magnetohydrodynamics (see~\bibref{Sart}).

The particular point of this strategy is the mixture between the mass equation and the momentum equation to get the control of non-linear diffusive terms with density-dependent coefficients through the BD entropy.

We first multiply~\eqref{1} by $\dsp \vph'(\vrh)=\frac{\mu'(\vrh)}{\vrh}$ to get:
\[
\pa_t\vph(\vrh)+\u\cdot\na\vph(\vrh)+\mu'(\vrh)\div (\u) = 0.
\]
Then we derive with respect to the space variables:
\[
\pa_t\na\vph(\vrh)+(\u\cdot\na)\na\vph(\vrh)+\na \u:\na\vph(\vrh)+\na (\mu'(\vrh)\div (\u)) = 0.
\]
Let's now multiply by $2\vrh$, then, noting $\U=2\na\vph(\vrh)$ and using \eqref{1}, we write:
\[
\pa_t(\vrh\U)+\div(\vrh \u\otimes \U)+2\na \u:\na\mu(\vrh)+2\vrh\na\big (\mu'(\vrh)\div (\u) \big )=0.
\]
Rewriting the last two terms, one has
\begin{eqnarray} 
2\na \u:\na\mu(\vrh) & = & 2\div\big (\mu(\vrh)\na\u\big )-2\mu(\vrh)\na(\div (\u))\nonumber\\
                     & = & 2\div\big (\mu(\vrh)D(\u)\big )+2\div\big (\mu(\vrh)A(\u)\big )-2\na\big (\mu(\vrh)\div(\u)\big )+2\na\mu(\vrh)\div(\u),\nonumber\\
2\vrh\na\big (\mu'(\vrh)\div (\u)\big ) & = & 2\na\big (\vrh\mu'(\vrh)\div(\u)\big )-2\na\mu(\vrh)\div(\u), \nonumber
\end{eqnarray}
where we recall that~$D$ is the symmetric part of the gradient, and where~$A$ is the skew symmetric part of the gradient.
Then, summing with the momentum equation \eqref{2}, we get
\begin{equation*}
\begin{aligned}
\partial_t (\rho (\u+\U)) + \div (\rho \u \otimes (\u+\U)) + \nabla p(\rho) 
& = \div (2\mu(\rho)A(\u)) - r_0 \rho |\u| \u \\
& \qquad + \na\big ((2\vrh\mu'(\vrh)-2\mu(\vrh)-\la(\vrh))\div (\u)\big ).
\end{aligned}
\end{equation*}
Notice that the assumption \eqref{5} is now necessary to get the following interesting form:
\begin{equation}\label{BD}
  \partial_t (\rho (\u+\U)) + \div (\rho \u \otimes (\u+\U)) + \nabla p(\rho) = \div (2\mu(\rho)A(\u)) - r_0 \rho |\u| \u.
\end{equation}

\SUBSECTION{2.2}{A priori estimates}\label{sec:estimates}

To control the boundary terms in various integrations by parts, we introduce a lift of the velocity field. Let~$\tu$ be a regular function such that
\begin{equation*}
  \u = \tu \ \text{on $\partial \Omega$, \quad}
  \tu\cdot \n = 0 \; \text{on $\partial \Omega$ \quad and \quad}
  \div (\tu) = 0 \; \text{on $\Omega$.}
\end{equation*}

\SSUBSECTION{2.2.1}{Energy}

The energy estimate comes from the multiplication of~\eqref{2} by $\u-\tu$. Using equation~\eqref{1} and the boundary conditions on~$\u$ and~$\tu$ we obtain
\begin{equation}\label{eq:energy}
\begin{aligned}
&  \frac{d}{dt} \int_\Omega \left( \rho \frac{|\u|^2}{2} + Q(\rho) \right)
+ \int_\Omega 2\mu(\rho) |D(\u)|^2 
+ \int_\Omega \lambda(\rho) |\div (\u)|^2
+ r_0 \int_\Omega \rho \, |\u|^{3} \\
&  \quad  = \frac{d}{dt} \int_\Omega \left( \rho \u\cdot \tu \right)
- \int_\Omega (\rho \u\otimes \u) : \nabla \tu
+ \int_\Omega 2\mu(\rho) D(\u):D(\tu) 
+ r_0 \int_\Omega \rho |\u| \u\cdot\tu,
\end{aligned}
\end{equation}
where $Q'(\rho):=\Pi(\rho)$ and $\rho \Pi'(\rho):=p(\rho)$.

\SSUBSECTION{2.2.2}{BD entropy}

The BD entropy estimate comes from the multiplication of~\eqref{BD} by $\u-\tu+\U$. 
Using equation~\eqref{1} and the boundary conditions on~$\u$, $\tu$ and $\rho$ we obtain
\begin{equation}\label{eq:BD}
\begin{aligned}
&  \frac{d}{dt} \int_\Omega \left( \rho \frac{|\u+\U|^2}{2} + Q(\rho) \right)
+ \int_\Omega 2\mu(\rho) |A(\u)|^2
+ \int_\Omega \nabla P(\rho) \cdot \nabla \fy(\rho) 
+ r_0 \int_\Omega \rho \, |\u|^{3} \\
&  \qquad  
= \frac{d}{dt} \int_\Omega \left( \rho (\u+\U)\cdot \tu \right)
- \int_\Omega (\rho \u\otimes (\u+\U)) : \nabla \tu
- \int_\Omega 2\mu(\rho) A(\u):A(\tu)\\
& \hspace{7cm}
+ r_0 \int_\Omega \rho |\u| \u\cdot\tu
- r_0 \int_\Omega \rho |\u| \u\cdot\U.
\end{aligned}
\end{equation}
One of the main interests of this estimate is that not only it makes it possible to have a control on $\U$ {\it via} the control of $\rho (\u+\U)$ but also that the ``pressure'' term $\nabla p(\rho) \cdot \nabla \fy(\rho)$ is very rich. \par
Separating the pressure into two terms : $p=p_h + p_c$, see assumption~\eqref{8}, we write
\[
\int_\Omega \nabla p \cdot \nabla \fy(\vrh) = \int_\Omega \nabla p_h(\rho) \cdot \nabla \fy(\vrh)+\int_\Omega \nabla p_c(\rho) \cdot \nabla \fy(\vrh) =: I_h + I_c.
\]
About the term $I_h$, we use the definition of $p_h(\vrh)=a\vrh^\gamma$ and find
\[
I_h = a\, \gamma \int_\Omega \mu'(\vrh) \vrh^{\gamma-2}|\nabla \vrh|^2.
\]

We write the term $I_c$ as follows
\begin{equation*}
I_c = 
\int_\Omega p_c'(\rho) \frac{\mu'(\rho)}{\rho} |\nabla \rho|^2 \, \mathds{1}_{\rho<\rho_*}
+ \int_\Omega p_c'(\rho) \frac{\mu'(\rho)}{\rho} |\nabla \rho|^2 \, \mathds{1}_{\rho>\rho_*}.
\end{equation*}
Using the assumptions~\eqref{9} and~\eqref{10} we obtain
\begin{equation*}
I_c \geq
\frac{1}{c_2}\int_\Omega \rho^{-\alpha-2} \mu'(\rho) |\nabla \rho|^2 \, \mathds{1}_{\rho<\rho_*}
- \frac{a\, \gamma}{2} \int_\Omega \rho^{\gamma-2} \mu'(\rho) |\nabla \rho|^2.
\end{equation*}
Moreover, for $\rho<A$ we can use the assumption~\eqref{6} on~$\mu'$ and deduce
\begin{equation*}
I_c \geq
\frac{c_0'}{c_2}\int_\Omega \rho^{n-\alpha-3} |\nabla \rho|^2 \, \mathds{1}_{\rho<\min(\rho_*,A)}
- \frac{a\, \gamma}{2} \int_\Omega \rho^{\gamma-2} \mu'(\rho) |\nabla \rho|^2.
\end{equation*}
Adding the two contributions $I_h$ and $I_c$ we obtain
\begin{equation*}
\int_\Omega \nabla p \cdot \nabla \fy(\vrh) 
\geq
\frac{c_0'}{c_2 \, M^2} \int_\Omega \Big |\nabla \left( \xi(\vrh)^M \right)\Big |^2 
+ \frac{a\, \gamma}{2} \int_\Omega \rho^{\gamma-2} \mu'(\rho) |\nabla \rho|^2.
\end{equation*}
where $M:=\frac{-\alpha-1+n}{2}$ and $\xi$ being taken such that $\xi(\vrh)=\vrh$ for $\vrh\leq \frac{1}{2}\min(\vrh_*,A)$ and $\xi(\vrh)=0$ for $\vrh\geq \min(\vrh_*,A)$.
\par
Note that the contribution $\frac{a\, \gamma}{2} \int_\Omega \rho^{\gamma-2} \mu'(\rho) |\nabla \rho|^2$ allows us to control positive power of~$\rho$.
In fact, using assumptions~\eqref{6} and~\eqref{7}, that is separating small and large densities, we show that
\begin{equation*}
\int_\Omega \rho^{\gamma-2} \mu'(\rho) |\nabla \rho|^2
\geq
C_0 \int_\Omega \vrh^{\gamma-3+n} |\nabla \vrh|^2,
\end{equation*}
where $C_0$ depends on~$c_0'$, $c_1'$ and~$A^{m-n}$.
We finaly obtain the following inequality.
\begin{equation}\label{14}
\int_\Omega \nabla p \cdot \nabla \fy(\vrh) \geq \frac{c_0'}{c_2 M^2}\int_\Omega \Big |\nabla \left( \xi(\vrh)^M \right)\Big |^2 + \frac{C_0\, a\,\gamma}{2N^2} \int_\Omega \Big |\nabla \left( \vrh^N \right) \Big |^2,
\end{equation}
where $M:=\frac{-\alpha-1+n}{2}<0$ and $N:=\frac{\gamma+n-1}{2}>0$.

\SSUBSECTION{2.2.3}{Control of integral terms}\label{non-stationary}

In this part, we are particularily interested in the case $r_0=0$. 
In this case, the evolution terms are enough to control all the other terms. Nevertheless, an additional friction term with $r_0>0$ naturally preserves the following calculations.\par
Integrating with respect to the time $t\in [0,T]$ the energy estimate~\eqref{eq:energy}, we obtain
\begin{equation}\label{eq:energyint}
\begin{aligned}
&  \int_\Omega \left( \rho(T) \frac{|\u(T)|^2}{2} + Q(\rho(T)) \right)
+ \int_0^T \int_\Omega 2\mu(\rho) |D(\u)|^2 
+ \int_0^T \int_\Omega \lambda(\rho) |\div (\u)|^2\\
& \hspace{2cm} = 
\int_\Omega \left( \frac{|\m_0|^2}{2\rho_0} + Q(\rho_0) \right)
+ \int_\Omega \rho(T) \u(T)\cdot \tu
- \int_\Omega  \m_0 \cdot \tu  \\
& \hspace{2cm} \qquad
- \int_0^T \int_\Omega (\rho \u\otimes \u) : \nabla \tu
+ \int_0^T \int_\Omega 2\mu(\rho) D(\u):D(\tu).
\end{aligned}
\end{equation}
Each term of the right hand side member is controlled as follow:
\begin{equation*}
\begin{aligned}
\hspace{-4cm}\bullet \;
\left| \int_\Omega \rho(T) \u(T)\cdot \tu \right|
& \leq 
  \frac{1}{2} \int_\Omega \rho(T) \frac{|\u(T)|^2}{2} + \int_\Omega \rho(T) |\tu|^2\\
& \leq 
  \frac{1}{2} \int_\Omega \rho(T) \frac{|\u(T)|^2}{2} + |\tu|_\infty^2 M_0
\end{aligned}
\end{equation*}
\begin{equation*}
\begin{aligned}
\hspace{-4.4cm}\bullet \;
  \left| \int_0^T \int_\Omega (\rho \u\otimes \u) : \nabla \tu \right|
& \leq
  |\nabla \tu|_\infty \int_0^T \int_\Omega \rho |\u|^2
\end{aligned}
\end{equation*}
\begin{equation*}
\begin{aligned}
\bullet \;
  \left| \int_0^T \int_\Omega 2\mu(\rho) D(\u):D(\tu) \right|
& \leq 
  \int_0^T \int_\Omega \mu(\rho) |D(\u)|^2  + \int_0^T \int_\Omega \mu(\rho) |D(\tu)|^2 \\
& \leq 
  \int_0^T \int_\Omega \mu(\rho) |D(\u)|^2  + |D(\tu)|_\infty^2 \int_0^T \int_\Omega \mu(\rho).
\end{aligned}
\end{equation*}
We deduce from~\eqref{eq:energyint} that
\begin{equation}\label{eq:energyint1}
\begin{aligned}
&  \int_\Omega \left( \rho(T) \frac{|\u(T)|^2}{4} + Q(\rho(T)) \right)
+ \int_0^T \int_\Omega \mu(\rho) |D(\u)|^2 
+ \int_0^T \int_\Omega \lambda(\rho) |\div (\u)|^2 \\
&  \hspace{4cm}  \leq
C(\m_0,\rho_0,\tu)
+ |\nabla \tu|_\infty \int_0^T \int_\Omega \rho |\u|^2
+ |D(\tu)|_\infty^2 \int_0^T \int_\Omega \mu(\rho).
\end{aligned}
\end{equation}
In the same way, integrating with respect to the time $t\in [0,T]$ the BD entropy estimate~\eqref{eq:BD}, we obtain (recall that in this subsection the friction coefficient is assume to be zero)
\begin{equation}\label{eq:BDint}
\begin{aligned}
&  \int_\Omega \left( \rho(T) \frac{|\u(T)+\U(T)|^2}{2} + Q(\rho(T)) \right)
+ \int_0^T \int_\Omega 2\mu(\rho) |A(\u)|^2
+ \int_0^T \int_\Omega \nabla p(\rho) \cdot \nabla \fy(\rho) \\
&  \hspace{1cm}  
= \int_\Omega \left( \rho_0 \frac{|\u_0+\U_0|^2}{2} + Q(\rho_0) \right)
+ \int_\Omega \left( \rho(T) (\u(T)+\U(T))\cdot \tu \right)
- \int_\Omega \left( \rho_0 (\u_0+\U_0)\cdot \tu \right)\\
& \hspace{4cm}
- \int_0^T \int_\Omega (\rho \u\otimes (\u+\U)) : \nabla \tu
- \int_0^T \int_\Omega 2\mu(\rho) A(\u):A(\tu).
\end{aligned}
\end{equation}
The right-hand side members are estimated in the same way that for obtaining the estimate~\eqref{eq:energyint1}.
We obtain
\begin{equation}\label{eq:BDint1}
\begin{aligned}
&  \int_\Omega \left( \rho(T) \frac{|\u(T)+\U(T)|^2}{4} + Q(\rho(T)) \right)
+ \int_0^T \int_\Omega \mu(\rho) |A(\u)|^2
+ \int_0^T \int_\Omega \nabla p(\rho) \cdot \nabla \fy(\rho) \\
&  \hspace{1cm}  
\leq C(\m_0,\rho_0,\tu)
+ \frac{|\nabla \tu|_\infty}{2} \int_0^T \int_\Omega \rho |\u|^2
+ \frac{|\nabla \tu|_\infty}{2} \int_0^T \int_\Omega \rho |\u+\U|^2
+ |A(\tu)|_\infty^2 \int_0^T \int_\Omega \mu(\rho).
\end{aligned}
\end{equation}
%

\SSUBSECTION{2.2.4}{Gronwall argument}\label{Gronwall}

Putting \eqref{14}, \eqref{eq:energyint1} and \eqref{eq:BDint1} together we get 
\begin{equation*}
\begin{aligned}
& \int_\Omega \left( \rho(T) \frac{|\u(T)|^2}{4} +\rho(T) \frac{|\u(T)+\U(T)|^2}{4}+ 2Q(\rho(T)) \right)
+ \int_0^T \int_\Omega \mu(\rho) |D(\u)|^2 
+ \int_0^T \int_\Omega \mu(\rho) |A(\u)|^2 \\
& \qquad + \int_0^T \int_\Omega \lambda(\rho) |\div (\u)|^2
+\frac{C_0}{c_2 M^2}\int_0^T\int_\Omega \Big |\nabla \left( \xi(\vrh)^M \right)\Big |^2 +
\frac{C_0\, a\,\gamma}{2N^2} \int_\Omega \Big |\nabla \left( \vrh^{N} \right) \Big |^2 \\
& \qquad \qquad \leq
C(\m_0,\rho_0,\tu)
+ C(|\na\tu|_\infty)\Big [\int_0^T \int_\Omega \rho |\u|^2
+ \int_0^T \int_\Omega \rho |\u+\U|^2
\Big ]+ |\na\tu|_\infty^2\int_0^T \int_\Omega \mu(\rho).
\end{aligned}
\end{equation*}
Since the gradient of both positive and negative powers of the density appear on the left hand side (recall that $M=\frac{n-\alpha-1}{2}<0$ and $N=\frac{\gamma+n-1}{2}>0$) and since the density is constant on~$\pa\Om$, we can insure, thanks to Poincar\'e, that $\int_0^T\int_\Om \mu(\rho)$ is controlled, {\it via} assumptions~\eqref{9} and~\eqref{10}, by the pressure terms of the left hand side. 
Then, using a Gronwall argument, we deduce that all the left hand side terms of this last inequality are bounded and we can write the following estimates:

\begin{eqnarray}
\|\rho\|_{L^\infty(0,T;L^\gamma(\Omega))} & \leq &  c(\m_0,\rho_0,\tilde{\u}),\nonumber\\
\|\sqrt{\rho}\u\|_{L^\infty(0,T;L^2(\Omega))} & \leq &  c(\m_0,\rho_0,\tilde{\u}),\nonumber\\
\|\sqrt{\rho} \nabla \fy(\rho)\|_{L^\infty(0,T;L^2(\Omega))} & \leq &  c(\m_0,\rho_0,\tilde{\u}),\nonumber\\
\|\sqrt{\mu(\rho)}\nabla \u\|_{L^2((0,T)\times \Omega)} & \leq &  c(\m_0,\rho_0,\tilde{\u}),\nonumber\\
\|\na\xi(\rho)^{\frac{n-\al-1}{2}}\|_{L^2(0,T;L^2(\Omega))} & \leq &  c(\m_0,\rho_0,\tilde{\u}),\nonumber\\
\|\na\rho^{\frac{\ga+n-1}{2}}\|_{L^2(0,T;L^2(\Omega))} & \leq &  c(\m_0,\rho_0,\tilde{\u}).\nonumber
\end{eqnarray}

\SUBSECTION{2.3}{Approximate solutions and compactness}\label{sub:regul}

The preceding a priori estimates are the key ingredient of the existence result. 
As soon as approximate solutions satisfy such estimates, compactness properties make it possible to extract a subsequence that converges to a weak solution of the initial model.
The compactness arguments are exactly those given in the periodic case or the whole space in \bibref{BD} and more recently in the bounded case, see~\cite{BDG}, that is why we will not detail it here.

Let us just say some words about the sequences of suitably smooth approximate solutions to the compressible Navier-Stokes equations \eqref{1}--\eqref{2} that preserve the estimates obtained in the previous section.

The construction scheme of approximate solutions, using on additional regularizing effects such as capillarity, is provided in~\cite{BD3}.
We thus introduce some modified Navier-Stokes equations for $(\vrh_{\al,\be},\u_{\al,\be})$, always denoted~$(\vrh,\u)$ for sake of simplicity, depending on the regularizing parameters $\al$ and $\be$:
\begin{eqnarray}
  \pa_t\vrh+\div (\vrh \u) & = & 0, \label{15} \\
  \pa_t(\vrh \u)+\div (\vrh \u\otimes \u) -\div\sigma- r_0 \rho |\u| \u
  -\beta \vrh \nabla (\mu'(\vrh)\lap ^s \mu(\vrh)) + \alpha \lap^2 \u & = & 0, \label{16}
\end{eqnarray}
where the conditions \eqref{5}--\eqref{10} are supposed to be satisfied.\par
These regularizations allow to use some classical result in order to prove the existence of smooth solutions. 
The remaining work consists in showing that the additional terms depending on~$\alpha$ and on~$\beta$ lead to some weak solutions of our initial model \eqref{1}--\eqref{2}.\par
Notice that we may not modify \eqref{15} because the BD entropy is very closely related to the mass equation and some regularizing term in \eqref{15} could cancel equation \eqref{BD}.
For this model, since energy and BD identities are preserved, the stability arguments given in~\bibref{BD3} and~\bibref{BD} lead to our existence result cited in the Theorem~\ref{th-non-stationary}.

\SECTION{3}{Sketch of proof of Theorem~\ref{th-stationary} (stationary case)}\label{sec-3}

The proof of Theorem~\ref{th-non-stationary} has been managed in the general case $r_0\geq 0$. 
The only two points that have to be cleared in the stationary context concern the BD structure and the control of integral terms in the energy and BD formula.
To control these terms we assume in this part, as announced in Theorem~\ref{th-stationary}, that $r_0>0$.
Moreover, in the three-dimensional case, we will assume the additive condition~\eqref{condition}.
\par\noindent
These two additive conditions (the condition on~$r_0$ and the condition~\eqref{condition}) will be used since in the stationary case we can not use Gronwall type arguments.

\SUBSECTION{3.1}{BD structure}

Bringing some modification to the mass equation could cancel the BD structure, that is why it is not clear that Subsection \subsectionref{2.1} can be directly adapted. 
For instance, it is dedicated to the failure for any semi-stationary model, whereas stationary conditions for both mass and momentum equation lead, following the same steps as for the equation~\eqref{BD}, to a similar equation:
\[
\div (\rho \u \otimes (\u+\U)) + \nabla p(\rho) = \div (2\mu(\rho)A(\u)) - r_0 \rho |\u| \u,
\]

\SUBSECTION{3.2}{Control of integral terms}\label{stationary}

In the stationary case, we will use the friction term to obtain a ``good'' estimate. More precisly in this case the energy estimate~\eqref{eq:energy} and the BD entropy estimate~\eqref{eq:BD} respectively write
\begin{equation}\label{eq:energy-station}
\begin{aligned}
& \int_\Omega 2\mu(\rho) |D(\u)|^2 
+ \int_\Omega \lambda(\rho) |\div (\u)|^2
+ r_0 \int_\Omega \rho \, |\u|^{3} \\
&  \hspace{2cm}  
= - \int_\Omega (\rho \u\otimes \u) : \nabla \tu
+ \int_\Omega 2\mu(\rho) D(\u):D(\tu) 
+ r_0 \int_\Omega \rho |\u| \u\cdot\tu,
\end{aligned}
\end{equation}
\begin{equation}\label{eq:BD-station}
\begin{aligned}
& \int_\Omega 2\mu(\rho) |A(\u)|^2
+ \int_\Omega \nabla P(\rho) \cdot \nabla \fy(\rho) 
+ r_0 \int_\Omega \rho \, |\u|^{3} \\
&  \hspace{2cm}  
= - \int_\Omega (\rho \u\otimes (\u+\U)) : \nabla \tu
- \int_\Omega 2\mu(\rho) A(\u):A(\tu)\\
& \hspace{5cm}
+ r_0 \int_\Omega \rho |\u| \u\cdot\tu
- r_0 \int_\Omega \rho |\u| \u\cdot\U.
\end{aligned}
\end{equation}
We estimate the terms of right-hand sides again (the nonhere detailed terms are exactly treated as in the nonstationary case).
The constant~$C$ which appears does not depend on physical constants such~$\Omega$, $r_0$, $\tu$...
\begin{equation*}
\begin{aligned}
~~~\bullet \;
\left| \int_\Omega (\rho \u\otimes \u) : \nabla \tu \right|
& \leq 
\frac{r_0}{4} \int_\Omega \rho \, |\u|^{3} + \frac{C}{r_0}|\nabla \tu|_\infty^{3} \int_\Omega \rho,
\end{aligned}
\end{equation*}
\begin{equation*}
\begin{aligned}
\bullet \;
\left| r_0 \int_\Omega \rho |\u| \u\cdot\tu \right|
& \leq 
\frac{r_0}{4} \int_\Omega \rho \, |\u|^{3} + C\, r_0|\tu|_\infty^{3} \int_\Omega \rho.
\end{aligned}
\end{equation*}
The only two terms which seem more difficult to control are the following
\begin{equation*}
T_1 = \int_\Omega (\rho \u\otimes \U) : \nabla \tu
\quad \text{and}\quad
T_2 = r_0 \int_\Omega \rho |\u| \u\cdot\U.
\end{equation*}
$\bullet \;$
Using the definition of~$\U$ and of~$\fy$, and using an integration by part, since $\div(\tu)=0$, we obtain
\begin{equation*}
T_1 =  \int_\Omega 2(\u \otimes \nabla \mu(\rho)):\nabla \tu = -\int_\Omega 2\mu(\rho) \, (\nabla \u)^T : \nabla\tu + \int_{\partial \Omega} 2\mu(\rho) \, (\u\cdot \nabla \tu )\cdot \n.
\end{equation*}
Since $\tu\cdot \n=0$ and $\u=\tu=\u^b$ on $\partial \Omega$ we write 
\begin{equation*}
(\u\cdot \nabla \tu )\cdot \n = \u_i(\partial_i \tu_j)\n_j = \u_i \partial_i (\tu_j \n_j) -  \u_i (\partial_i \n_j) \tu_j = -\u^b\cdot \nabla \n \cdot \u^b = \mathrm{II}(\u^b),
\end{equation*}
where $\mathrm{II}$ is the second fundamental form of $\partial \Omega$. By definition of $A(\u)$ and $D(\u)$ we have
\begin{equation*}
  (\nabla \u)^T : \nabla\tu
  = (D(\u)-A(\u)):(D(\tu)+A(\tu))
  = D(\u):D(\tu) - A(\u):A(\tu).
\end{equation*}
Hence we obtain
\begin{equation*}
  T_1 \leq 
  \int_\Omega \mu(\rho) |D(\u)|^2 + |D(\tu)|_\infty^2 \int_\Omega \mu(\rho)
  + \int_\Omega \mu(\rho) |A(\u)|^2 + |A(\tu)|_\infty^2 \int_\Omega \mu(\rho)
  + 2\mu(\rho)^b \mathrm{II}(\u^b).
\end{equation*}
Finally, using the assumptions~\eqref{6} and~\eqref{7} we can control~$\mu$ with~$\rho$ as follows
\begin{equation*}
  \int_\Omega \mu(\rho)
= \int_\Omega \mu(\rho)\, \mathds{1}_{\rho<A}
+ \int_\Omega \mu(\rho)\, \mathds{1}_{\rho>A}
\leq \frac{1}{c_0} \int_\Omega \rho^n\, \mathds{1}_{\rho<A}
+ \frac{1}{c_1} \int_\Omega \rho^m\, \mathds{1}_{\rho>A}
\leq \frac{|\Omega|A^{n}}{c_0}
+ \frac{1}{c_1} \int_\Omega \rho^m.
\end{equation*}
We obtain
\begin{equation*}
  T_1
\leq \int_\Omega \mu(\rho) |D(\u)|^2 
+ \int_\Omega \mu(\rho) |A(\u)|^2 
+ \frac{|\Omega|\, A^n \, |\nabla \tu|_\infty^2}{c_0}
+ \frac{|\nabla \tu|_\infty^2}{c_1} \int_\Omega \rho^m 
+ 2\mu(\rho)^b \mathrm{II}(\u^b).
\end{equation*}
$\bullet \;$
For the term~$T_2$, since $\rho \U=2\nabla \mu(\rho)$ and $\u\cdot \n=0$ on $\partial \Omega$, we obtain by integration by part
\begin{equation*}
  T_2 = -\int_\Omega 2r_0\, \mu(\rho) \left( |\u|\div(\u) + \frac{\u}{|\u|}\cdot (\u\cdot \nabla)\u \right)
\leq 4 r_0 \int_\Omega \mu(\rho) |\u| |\nabla \u|.
\end{equation*}
By the Young inequality, we obtain
\begin{equation*}
  T_2 \leq \int_\Omega \mu(\rho) |\nabla \u|^2 + 4r_0^2 \int_\Omega \mu(\rho) |\u|^2.
\end{equation*}
Using assumptions~\eqref{6} and~\eqref{7}, the fact that $n\geq 2/3$ and the Young inequality, we succesively deduce that
\begin{equation*}
\begin{aligned}
4r_0^2 \int_\Omega \mu(\rho) |\u|^2 
& \leq 4r_0^2 \int_\Omega \mu(\rho) |\u|^2 \, \mathds{1}_{\rho<A} 
+ 4r_0^2 \int_\Omega \mu(\rho) |\u|^2 \, \mathds{1}_{\rho>A} \\
& \leq \frac{4r_0^2\, A^{n-2/3}}{c_0} \int_\Omega \rho^{2/3} |\u|^2
+ \frac{4r_0^2}{c_1} \int_\Omega \rho^m |\u|^2 \\
& \leq \frac{r_0}{4} \int_\Omega \rho |\u|^3 
+ \frac{C\, r_0^4\, A^{3n-2}}{c_0^3}
+ \frac{C\, r_0^4}{c_1^3} \int_\Omega \rho^{3m-2}.
\end{aligned}
\end{equation*}
Consequently we majore~$T_2$ as follows
\begin{equation*}
T_2
\leq \int_\Omega \mu(\rho) |\nabla \u|^2 
+ \frac{r_0}{4} \int_\Omega \rho |\u|^3 
+ \frac{C\, r_0^4\, A^{3n-2}}{c_0^3}
+ \frac{C\, r_0^4}{c_1^3} \int_\Omega \rho^{3m-2}.
\end{equation*}
With the preceding estimates, the sum of the equalities~\eqref{eq:energy-station} and~\eqref{eq:BD-station} is written
\begin{equation}\label{estimate}
\begin{aligned}
& \int_\Omega \mu(\rho) |D(\u)|^2 
+ \int_\Omega \mu(\rho) |A(\u)|^2
+ \int_\Omega \lambda(\rho) |\div (\u)|^2 \\
& \qquad +\frac{C_0}{c_2 M^2} \int_\Omega \Big |\nabla \left( \xi(\vrh)^M \right)\Big |^2 
+ \frac{C_0 a \gamma}{2} \int_\Omega \vrh^{\gamma-3+n} |\nabla \vrh|^2 
+ r_0 \int_\Omega \rho |\u|^3 \\
& \hspace{2cm} \leq 
\Big( \frac{C}{r_0}|\nabla \tu|_\infty^3 + Cr_0 |\tu|_\infty^3 \Big) \int_\Omega \rho
+ \frac{|\nabla \tu|_\infty^2}{c_1} \int_\Omega \rho^m 
+ \frac{C\, r_0^4}{c_1^3} \int_\Omega \rho^{3m-2}
+ \text{Cte},
\end{aligned}
\end{equation}
with $\dsp \text{Cte}=\frac{|\Omega|\, |\nabla \tu|_\infty^2}{c_0} + \frac{C\, r_0^4\, A^{3n-2}}{c_0^3} + 2\mu(\rho)^b \mathrm{II}(\u^b)$.
\par\vspace{0.5cm}
We conclude this section by showing that all the terms of right-hand side of the equation~\eqref{estimate} (except the constant~$\text{Cte}$) can be controlled by the terms of the left-hand side.
This result is due to the control of the density {\it via} the term $\int_\Omega \vrh^{\gamma-3+n} |\nabla \vrh|^2 = \frac{1}{N^2}\int_\Omega \big |\nabla \left( \vrh^N \right) \big |^2$ where $N=\frac{\gamma+n-1}{2}$.
This term make it possible (using the Poincaré inequality) to control $\rho^{N}$ in~$H^1(\Omega)$.
>From the Sobolev embeddings, we deduce a control of~$\rho$ in~$L^{qN}(\Omega)$ (for all $q<+\infty$ in the $2$-dimensional case, and for all $q\leq \frac{2d}{d-2}$ in the $d$-dimensional case, $d>2$).
\par
If we assume that
\begin{equation}\tag{C1}\label{C1}
3m-2 < qN
\end{equation}
then all the integrals of the right-hand side of the equation~\eqref{estimate} are controlled (since $m\geq 1$, that is $3m-2\geq m\geq 1$).
For instance, using the Young inequality, we write
\begin{equation*}
\frac{C\, r_0^4}{c_1^3} \int_\Omega \rho^{3m-2} \leq \delta_1 \int_\Omega \rho^{qN} + \delta_2,
\end{equation*}
where we can adapt the constant~$\delta_1$ such that the term $\delta_1 \int_\Omega \rho^{qN}$ is controlled by $\frac{C_0\, a\,\gamma}{2N^2} \int_\Omega \big |\nabla \left( \vrh^{N} \right) \big |^2$.

\SUBSECTION{3.4.0}{Stability of weak solutions}

In the stationary case, the lack of estimates implies that the stability of weak solutions is conditionned by some specific profiles for viscosities and pressure. 
Some relations between the corresponding coefficients $m, n, \al$ and $\ga$ may be considered. 
\par
Let us consider a sequence of weak solutions $\vrh_k,\u_k$ of the stationary equations \eqref{1stat}--\eqref{2stat}. 

\SSUBSECTION{3.4.1}{Estimates}

The preceding subsection leads to the following a priori estimates:

\begin{eqnarray}
\|\sqrt{\mu(\rho_k)}\nabla \u_k\|_{L^2(\Omega)} & \leq &  c(\Om,\tilde{\u}), \label{S1}\\
\|\sqrt{\la(\rho_k)}\div (\u_k)\|_{L^2(\Omega)} & \leq &  c(\Om,\tilde{\u}), \label{S11}\\
\|\na \left( \xi(\rho_k)^{M} \right) \|_{L^2(\Omega)} & \leq &  c(\Om,\tilde{\u}), \label{S2}\\ 
\|\na \left( \rho_k^{N} \right) \|_{L^2(\Omega)} & \leq &  c(\Om,\tilde{\u}), \label{S22}\\
\|\rho_k\u_k^3\|_{L^1(\Omega)} & \leq &  c(\Om,\tilde{\u}), \label{S3}
\end{eqnarray}
where $M=\frac{n-\alpha-1}{2}<0$ and $N=\frac{n+\gamma-1}{2}>0$.
\par
We are going to show that these estimates together with some compactness arguments lead to conclude that $(\vrh_k,\u_k)$ weakly converges to a solution $(\vrh,\u)$ of the system \eqref{1stat}--\eqref{2stat}.

\SSUBSECTION{3.4.2}{Compactnesses}

In order to cover the general case $d\in\{2,3\}$, we will keep a coefficient $q$ such that $H^1(\Om)\subset L^q(\Om)$ with continuous injection.
In the $d$-dimensional case (with $d>2$) we can choose any~$q$ such that $q\leq\nfrac{2d}{(d-2)}$ whereas in the $2$-dimensional case we can choose any~$q$ such that $q<+\infty$.
In the sequel, we will denote by $q$ such a real.
\par\vspace{0.3cm}
$\bullet$ {\bf Compacity on the density -}
The estimate~\eqref{S22} shows that the sequence $\vrh_k^N$ is bounded in $H^1(\Om)$.
Under the condition~\eqref{C1} and the fact that $3m-2\geq 1$ for all $m\geq 1$, we have $qN\geq 1$.
Consequently we obtain
\begin{equation}\label{S4}
  \vrh_k\rightarrow\vrh\ \ {\rm in}\ L^{qN}(\Om).
\end{equation}
In the same way, the estimate~\eqref{S2} shows that the sequence~$\vrh_k^M$ is bounded in $H^1(\Om)$ (recall that by definition we have~$M<0$).
We obtain
\begin{equation}\label{S5}
  \frac{1}{\vrh_k} \rightarrow \frac{1}{\vrh}\ \ {\rm in}\ L^{-qM}(\Om).
\end{equation}
We will note that $-qM\geq 1$ is satified in the $2$-dimensional case taking~$q$ large enough and in the $3$-dimensional case taking $q=6$ and using the assumptions given on page~\pageref{5} for~$\alpha$ and~$n$.
\par
By the conditions \eqref{6} and \eqref{7}, we obtain
\begin{equation}\label{S6}
  \sqrt{\mu(\vrh_k)}\rightarrow\sqrt{\mu(\vrh)}\ \ {\rm in}\ L^{\frac{2qN}{m}}(\Om),
\end{equation}
\begin{equation}\label{S7}
  \frac{1}{\sqrt{\mu(\vrh_k)}}\rightarrow\frac{1}{\sqrt{\mu(\vrh)}}\ \ {\rm in}\ L^{\frac{-2qM}{n}}(\Om).
\end{equation}
We will note that $\frac{2qN}{m}\geq 2$ (using the condition~\eqref{C1} and the fact that $3m-2 \geq m$ for all $m\geq 1$) and $\frac{-2qM}{n}\geq 2$ (since we previously prove that $-qM\geq 1$ and since $n<1$).
\par\vspace{0.3cm}
$\bullet$ {\bf Compacity on the velocity -}
On another hand, we know by \eqref{S1} that $\sqrt{\mu(\vrh_k)}\na \u_k$ is bounded in~$L^2(\Om)$ and thus weakly converges in $L^2(\Om)$.
>From the identity
\[
\na \u_k=\frac{1}{\sqrt{\mu(\vrh_k)}}\sqrt{\mu(\vrh_k)}\na\u_k,
\]
we also conclude that $\na\u_k$ is bounded in $L^r(\Om)$ with $\frac{1}{r}=\frac{1}{2}-\frac{n}{2qM}$.
We note that $r \geq 1$ since we have previously proved that $-qM \geq 1 > n$.
Moreover, since $M<0$, we also have $r<2$.
\par
Using Poincar\'e inequality we obtain a bound for the sequence~$\u_k$ in~$W^{1,r}_w(\Om)$.
Thanks to the compactness $W^{1,r}(\Om)\subset L^s(\Om)$ for $s<{\frac{rd}{d-r}}$, we obtain
\begin{equation}\label{S10}
  \u_k\rightarrow \u\ {\rm in}\ L^s(\Om),\ \forall s<{\frac{rd}{d-r}}.
\end{equation}
%

\SSUBSECTION{3.4.3}{Limit}

In this subsection, we show that we can pass to the limit when $k$ tend to $+\infty$ for the nonlinear term in the equation~\eqref{1stat} and~\eqref{2stat}.
The ``more nonlinear'' terms in these equations are the following ones:
\begin{equation*}
  T_1 = P(\rho_k),
 \quad
  T_2 = \rho_k |\u_k| \u_k
 \quad \text{and} \quad
  T_3 = \div \big (\mu(\vrh_k)\na \u_k\big ).
\end{equation*}
More precisely, the other nonlinear terms are $\div (\rho_k \u_k)$ and $\div(\rho_k \u_k\otimes \u_k)$ which convergences (in the sense of distributions on~$\Om$) are consequences of the convergence of~$T_2$, and $\na \big (\la(\vrh_k)\div(\u_k)\big )$ which convergence is similar to the convergence of~$T_3$.
\par\vspace{0.3cm}
$\bullet$ {\bf Convergence of the pressure term $T_1$ -}
Recall (see assumption~\eqref{8}) that the pressure is a sum of two pressures $p_h+p_c$.
\par
\begin{itemize}
\item[$\star$]
Since $p_h(\vrh_k)=a\vrh_k^\gamma$ the convergence of $\na p_h(\vrh_k)$ to $\na p_h(\vrh)$ in the sense of distributions on~$\Om$ comes from to convergence~\eqref{S4}:
\begin{equation}\label{S8}
\vrh_k^\ga\rightarrow \vrh^\ga\ {\rm in}\ L^{\frac{qN}{\ga}}(\Om).
\end{equation}
We will note that $\frac{qN}{\ga}\geq 1$.
More precisely, this condition is satisfied in the $2$-dimensional case (taking~$q$ large enough) and in the $3$-dimensional case taking $q=6$, and $n$ and $\gamma$ satisfying the assumptions given page~\pageref{5}.
\item[$\star$]
Then, we are interested in the convergence of the cold pressure term which writes as
\begin{equation*}
\na p_c(\vrh_k)
= \vrh_k^{-M-\al} \Big( \vrh_k^{M+\al}p'_c(\vrh_k)\na\vrh_k\mathds{1}_{\{\vrh_k\leq\vrh*\}}\Big)
+ \vrh_k^{\max\{\beta,\gamma\}-N} \Big( \vrh_k^{N-\max\{\beta,\gamma\}} p'_c(\vrh_k)\na\vrh_k\mathds{1}_{\{\vrh_k>\vrh*\}} \Big).
\end{equation*}
The only thing we have to obtain on the gradient of the cold pressure $\na p_c(\vrh_k)$ is its boundedness in some $L^t(\Om)$ space with $t\geq1$.
Recalling the assumptions~\eqref{9} and~\eqref{10} on the cold pressure, we know that $\vrh_k^{M+\al}|p'_c(\vrh_k)\na\vrh_k|\mathds{1}_{\{\vrh_k\leq\vrh*\}}$ and $\vrh_k^{N-\max\{\beta,\gamma\}} |p'_c(\vrh_k)\na\vrh_k|\mathds{1}_{\{\vrh_k>\vrh*\}}$ are bounded in $L^2(\Om)$, respectively by \eqref{S2} and \eqref{S22}. 
As a consequence, we can insure that $\na p_c(\vrh_k)$ is bounded in $L^t(\Om)$ with $t\geq 1$ as soon as $\vrh_k^{-M-\al}$ and $\vrh_k^{\max\{\beta,\gamma\}-N}$ are bounded in $L^2(\Om)$.
\par
Since $-M-\al<0$, using the convergence result~\eqref{S5} we get the expected information ``$\vrh_k^{-M-\al}$ is bounded in $L^2(\Om)$'' if we have
\begin{equation*}
2(\alpha+M)\leq -qM.
\end{equation*}
We note that this condition is satisfied in the $2$-dimensional case taking $q$ large enough and in the $3$-dimensional case taking $q=6$.
\par
In the same way, since $\max\{\beta,\gamma\}-N>0$, using the convergence result~\eqref{S4} we get the expected information ``$\vrh_k^{\max\{\beta,\gamma\}-N}$ is bounded in $L^2(\Om)$'' if we have
\begin{equation*}
2(\max\{\beta,\gamma\}-N)\leq qN.
\end{equation*}
As previously, we note that the condition $2 (\gamma - N) \leq qN$ is satisfied as well in the $2$-dimensional case as in the $3$-dimensional case.
Hence, we need the following condition $2 (\beta - N) \leq qN$ which can be written
\begin{equation}\tag{C2}\label{C2}
4 \beta \leq (q+2)(n+\gamma-1).
\end{equation}

\end{itemize}
\par\vspace{0.3cm}
$\bullet$ {\bf Convergence of the friction term $T_2$ -}
We write~$T_2 = \big( \rho_k |\u_k|^3 \big)^{\frac{2}{3}} \, \rho_k^{\frac{1}{3}}$.
Using the compacity (that is the strong convergence~\eqref{S4}) and the bound (see estimate~\eqref{S3}) on $\rho_k |\u_k|^3$ in $L^1(\Omega)$, we obtain
\begin{equation*}
  \big( \rho_k |\u_k|^3 \big)^{\frac{2}{3}} \, \rho_k^{\frac{1}{3}} \rightharpoonup f \, \rho^{\frac{1}{3}} \ \ {\rm in}\ L^{\frac{3}{2}}(\Om)\times L^{3qN}(\Om) \subset L^1(\Om),
\end{equation*}
where $f$ is the weak limit of $\big( \rho_k |\u_k|^3 \big)^{\frac{2}{3}}$ in $L^{\frac{3}{2}}(\Omega)$.
The last inclusion holds since $qN\geq 1$ (see condition~\eqref{C1}).
\par
To identify the limit~$f$, we  use the strong convergence for the density and the velocity:
\begin{equation*}
  \rho_k^{\frac{2}{3}} \to \rho^{\frac{2}{3}} \ \ {\rm in}\ L^{\frac{3qN}{2}}(\Om)
\qquad \text{and} \qquad
  |\u_k|^2 \to |\u|^2 \ \ {\rm in}\ L^s(\Om),\ \forall s<{\frac{rd}{2(d-r)}}.
\end{equation*}
We deduce that $f=\rho^{\frac{2}{3}}|\u|^2$ if $\frac{2}{3qN}+\frac{2(d-r)}{rd} < 1$.
We can show that this condition is satisfy in the two-dimensional case taking $q$ large enough and in the $3$-dimensional case taking $q=6$ (and using the fact that $N>\frac{1}{3}$, $M<\frac{-1}{2}$ and $n<1$).
\par
Consequently, the friction term $T_2$ satisfies
\begin{equation*}
  T_2 
= \rho_k |\u_k|^2 
= \big( \rho_k |\u_k|^3 \big)^{\frac{2}{3}} \, \rho_k^{\frac{1}{3}} 
\rightharpoonup  \rho^{\frac{2}{3}}|\u|^2\, \rho^{\frac{1}{3}} 
= \rho |\u|^2 \ \ {\rm in}\  L^1(\Om).
\end{equation*}
\par\vspace{0.3cm}
$\bullet$ {\bf Convergence of the viscous term $T_3$ -}
Through \eqref{S1} we obtain
\begin{equation*}
  \sqrt{\mu(\vrh_k)}\na \u_k\rightharpoonup g \ \ {\rm in}\ L^2(\Om).
\end{equation*}
To identify the limit~$g$, we use the strong convergence results~\eqref{S6} and~\eqref{S10}.
We get $g=\sqrt{\mu(\rho)} \nabla \u$ if we have the following condition
\begin{equation}
\label{C'3}
\frac{m}{N}-\frac{n}{M} \leq q.
\end{equation}
This condition is satisfied in the $2$-dimensional case taking $q$ large enough.
In the $3$-dimensional case, taking $q=6$ the condition~\eqref{C'3} is written (recall that $M=\frac{n-\al-1}{2}$ and $N=\frac{n+\gamma-1}{2}$)
\begin{equation*}
m<(\ga+n-1)\Big(3-\frac{n}{1+\al-n}\Big).
\end{equation*}
Since $n<1$ and $\alpha\geq 1$, we have $3-\frac{n}{1+\al-n} > 1$.
We deduce that the condition~\eqref{C'3} is contained in the condition~\eqref{C1} in the $3$-dimensional case.
\par
The viscous term $\mu(\vrh_k)\na \u_k$ is written $\sqrt{\mu(\vrh_k)} \big( \sqrt{\mu(\vrh_k)}\na \u_k \big)$ which converges in $L^1(\Omega)$ if $m \leq qN$, condition which is a consequence of the condition~\eqref{C1} since $3m-2\geq m$.
We obtain
\begin{equation*}
  \div \big (\mu(\vrh_k)\na \u_k\big )\rightarrow \div\big (\mu(\vrh)\na \u\big )\ \ {\rm in}\ \cD'(\Om).
\end{equation*}
%

\SUBSECTION{3.4}{Pressure and viscosity conditions}

Let's recapitulate all the conditions we need to get the integrabilities and compacities cited in the preceding subsections.
Recall that (see assumption on page~\pageref{5})
\begin{equation}\label{conditions-coeff0}
\ga\geq 1,\quad
\al\geq 1, \quad
m>1 \quad 
\text{and}\quad 
\frac{2}{3}<n<1.
\end{equation}
The additional conditions~\eqref{C1} and \eqref{C2} write as
\begin{equation*}
3m-2 < qN,
    \quad
4\beta \leq (q+2)(n+\gamma-1)
.
\end{equation*}
\vskip 5pt

{\bf In the two dimensional case}, $q$ can be chosen as large as we need, thus many inequalities are satisfied and Theorem \ref{th-stationary} holds only with conditions \eqref{conditions-coeff0}.

\vskip 5pt

{\bf In the three dimensional case}, $q$ is any number smaller than $6$.
Thus, Theorem \ref{th-stationary} holds with the following additive conditions on the pressure and viscosities coefficients:
\begin{equation*}
m < \ga+n-\frac{1}{3},
\quad
\beta \leq 2(\ga +n-1)
.
\end{equation*}
These conditions exactly correspond to the condition~\eqref{condition}.

\SECTION{4}{Behaviour in thin domains}\label{sec:thin}

In this part, we derive the compressible Reynolds equation.
Formally, this equation comes from the compressible Navier-Stokes equation in a thin domain, that is when one of the length is assumed to be smaller than the other directions.
The main applications of this kind of behavior relate to the field of lubrication (see the Introduction).
Within such a framework, the thin domain is of the following form
\[
\Omega_\eps = \{(\x,z)\in \R^{d-1}\times \R ~; \quad \x\in {\cal O}\subset \R^{d-1} \quad \text{and} \quad 0<z<\eps h(\x)\},
\]
where ${\cal O}$ is a bounded domain in $\R^{d-1}$ and the height $h:{\cal O} \rightarrow \R$ is a regular and periodic function.
Note that to be able to define a periodical function, the domain $\mathcal O$ must be rectangular.
In the case of the dimension $d=2$ this is not a resctriction.
In the case of the upper dimension, this situation corresponds to realistic physical situations.
Moreover, it is possible to consider other conditions on the lateral boundaries. For all these aspects, consult thesis of S.~Martin~\cite{M}, as well as named references.\par
We assume that $h \geq h_{\text{min}}>0$ and up to a normalization, we can assume that $h_{\text{min}}=1$.
The size of the bounded domain ${\cal O} \subset \R ^{d-1}$ is assumed to be of order~$1$.
The non-dimensional number~$\eps$ corresponds to the characteristic ratio between the characteristic lenghts of $\cal O$ and the characteristic height~$\eps h$.
%
\begin{figure}[htbp]
\begin{center}
{\psfrag{Omega}{$\Omega_\eps$}\psfrag{x}{$\x\in \R^{d-1}$}\psfrag{z}{$z\in \R$}\psfrag{e}{$\eps$}\psfrag{rhobas}{$\u=(V,0)$ \qquad $\rho=\rho_b$}\psfrag{rhohaut}{$\u=(0,0)$ \qquad $\rho=\rho_t$}\psfrag{O1}{$\mathcal O (1)$}\psfrag{Oe}{$\mathcal O (\eps)$}
\includegraphics[width=12cm]{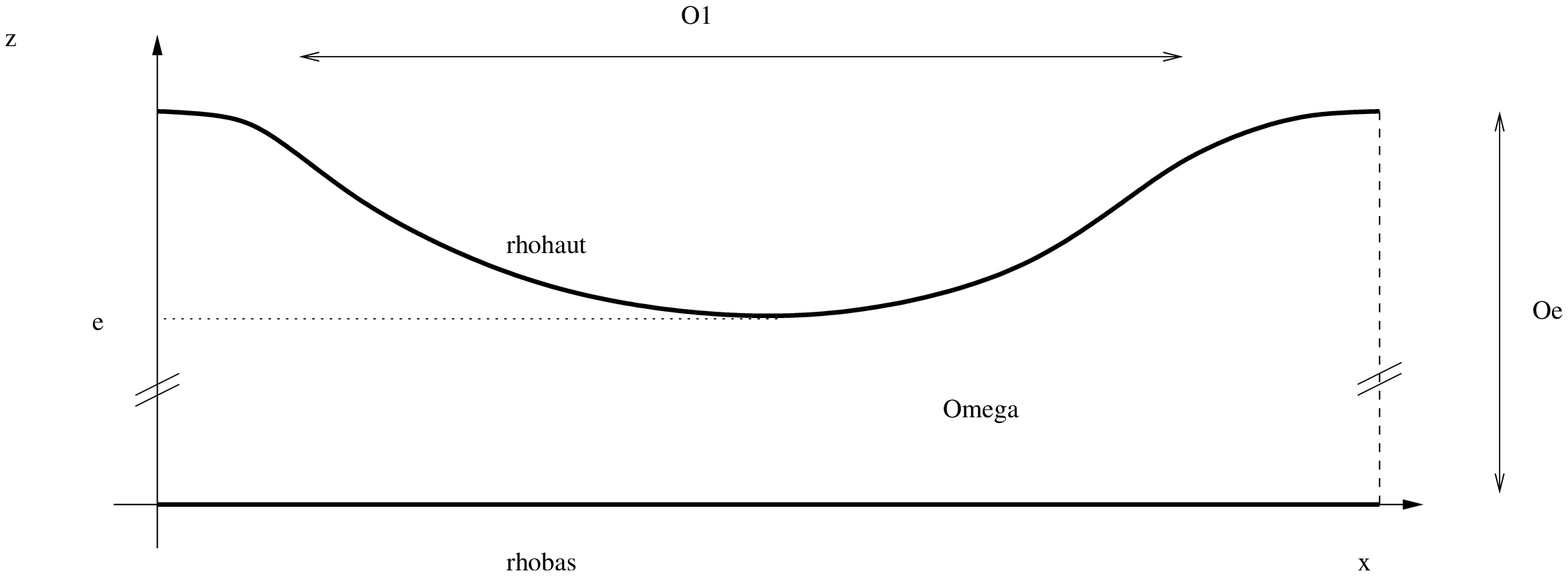}
}
\label{fig-domain}
\end{center}
\end{figure}
%

{\bf Boundary conditions on $\Om_\eps$:}
According to the results of the preceding parts (and according to the periodic results, see for instance~\cite{BD}), the boundary conditions which we impose are the following
\begin{enumerate}
\item[(i) -] Periodic conditions for the velocity and the density on the lateral boundaries (i.e. for $\x\in \partial \mathcal O$).
\item[(ii) -] Dirichlet conditions for the velocity on the top and bottom surface
\[\u=(V,0)\in \R^{d-1}\times \R \quad \text{for $z=0$}, \qquad \u=(0,0) \quad \text{for $z=h(\x)$}.\]
\item[(iii) -] Constant density on each connex component
\[\rho=\rho_b\in \R \quad \text{for $z=0$}, \qquad \rho=\rho_t\in \R \quad \text{for $z=h(\x)$}.\]
\end{enumerate}
The goal of this part is thus to justify in a rigorous way the compressible Reynolds equations, i.e. to determine the limit when~$\eps$ tends to~$0$ of the stationary compressible Navier-Stokes equations~\eqref{1stat}--\eqref{2stat}.
\REMARK
The result which is shown here concerns the justification of the Reynolds equation from the stationary Navier-Stokes compressible equations. Of course, the same method would allow to give a rigorous justification of the non-stationary Reynolds equation from the non-stationary Navier-Stokes equations.
\ENDREMARK

\SUBSECTION{4.1}{Rescaled equations}

In such a domain the unknowns of the equations~\eqref{1stat}--\eqref{2stat}, i.e. velocity and density, depending on~$\eps$ are denoted with a subscript~$\u_\eps$ and~$\rho_\eps$.
The first stage consists in rewriting these equations~\eqref{1stat}--\eqref{2stat} in a domain independent of~$\eps$. For that, we introduce the change of variable $Z=z/\eps$.
We define the rescaled domain
\[
\Omega = \{(\x,Z)\in \R^{d-1}\times \R ~; \quad \x\in {\cal O}\subset \R^{d-1} \quad \text{and} \quad 0<Z<h(\x)\}.
\]
In worries of simplifications, computations and notations used later will be made in dimension $d=2$. 
In this case, the impose velocity~$V$ is a real number, which will be assumed to be positive: $V>0$.
If the three dimensional case ($d=3$) is really different (for instance when we use the classical Sobolev injections) we shall apparently refer to it.\par
In the studied context (for example that of lubrication), we know that the pressure depends on the thickness~$\eps$ of the domain as $1/\eps^2$ (see~\cite{M} and the cited references). We define a normalized pressure~$P_\eps$ by $P_\eps=\eps^2 p_\eps$.
In the same way, if horizontal velocity is of order~$1$ (this order of magnitude depends in fact on the size of the velocity imposed on the boundaries of the domain, here we suppose that $V$ is of order~$1$) then vertical velocity will be of order $\eps$: $\u_\eps=(v_\eps,\eps w_\eps)$.
We can obtain the following equations in~$\Omega$:
\begin{eqnarray}
  \pa_x (\vrh_\vep v_\vep) + \pa_Z (\vrh_\vep w_\vep) & = & 0,\label{18}\\
  \vrh_\vep v_\vep\pa_x v_\vep + \vrh_\vep w_\vep\pa_Z v_\vep & = & 2\pa_x(\mu(\vrh_\vep) \pa_x v_\vep) + \frac{1}{\vep^2}\pa_Z\big (\mu(\vrh_\vep) \pa_Z v_\vep\big )+ \pa_Z\big (\mu(\vrh_\vep) \pa_x w_\vep\big ) \label{19}\\
   & & + \pa_x\big (\lambda(\vrh_\vep) (\pa_x v_\vep + \pa_Z w_\vep)\big ) - \frac{1}{\vep^2} \pa_x P(\vrh_\vep)-r_0\vrh_\vep (v_\vep^2+\vep^2 w_\vep^2)^{\frac{1}{2}} v_\vep, \nonumber\\
  \vrh_\vep v_\vep\pa_x(\vep w_\vep) +\vrh_\vep w_\vep\pa_Z (\vep w_\vep) & = &  \frac{1}{\vep}\pa_x\big (\mu(\vrh_\vep)\pa_Z v_\vep\big )+ \vep\pa_x\big (\mu(\vrh_\vep) \pa_x w_\vep\big ) + \frac{2}{\vep}\pa_Z\big (\mu(\vrh_\vep) \pa_Z w_\vep\big )\label{20}\\
     & & + \frac{1}{\vep} \pa_Z\big (\lambda(\vrh_\vep) (\pa_x v_\vep + \pa_Z w_\vep)\big )-\frac{1}{\vep^3} \pa_Z P(\vrh_\vep)-\vep r_0\vrh_\vep (v_\vep^2+\vep^2 w_\vep^2)^{\frac{1}{2}} w_\vep. \nonumber
\end{eqnarray}
As we have seen it in the Subsection \subsectionref{2.1}, we are also interested in other particular forms of these equations. 
Refering to \eqref{BD} and noting $\U_\vep=(V_\vep,W_\vep)=(\frac{2\pa_x\mu(\vrh_\vep)}{\vrh_\vep},\frac{2\pa_Z\mu(\vrh_\vep)}{\vep\vrh_\vep})$, we can write 
\begin{eqnarray}
  \vrh_\vep v_\vep\pa_x(v_\vep+V_\vep) + \vrh_\vep w_\vep\pa_Z(v_\vep+V_\vep) & = & \pa_Z(\mu(\vrh_\vep) \pa_x w_\vep)-\frac{1}{\vep^2}\pa_Z(\mu(\vrh_\vep) \pa_Z v_\vep) \nonumber \\
  & & -\frac{1}{\vep^2}\pa_x P(\vrh_\vep)-r_0\vrh_\vep (v_\vep^2+\vep^2 w_\vep^2)^{\frac{1}{2}} v_\vep, \label{BD-eps1}\\
 \vrh_\vep v_\vep \pa_x(\vep w_\vep+W_\vep) + \vrh_\vep w_\vep \pa_Z(\vep w_\vep+W_\vep) & = & \vep\pa_x(\mu(\vrh_\vep) \pa_x w_\vep)-\frac{1}{\vep}\pa_x(\mu(\vrh_\vep) \pa_Z v_\vep) \nonumber \\
  & & -\frac{1}{\vep^3}\pa_Z P(\vrh_\vep)-\vep r_0\vrh_\vep (v_\vep^2+\vep^2 w_\vep^2)^{\frac{1}{2}} w_\vep. \label{BD-eps2}
\end{eqnarray}

\SUBSECTION{4.2}{Convergence of the compressible Navier-Stokes equations to the compressible Reynolds equations}

\THEOREM{4}\label{th-convergence}
A solution $(\vrh_\vep,v_\vep,\vep w_\eps)$ of system \eqref{18}--\eqref{20} with conditions \eqref{5}--\eqref{condition+} satisfying the preceding boundary conditions (i), (ii) and (iii) converges to $(\vrh,v,w)$ in $L^{r_1}(\Om)\times \big (L_w^{r_2}(\Om)\big )^2$, for some $r_1,r_2>1$, when $\vep$ tends to~$0$.

At the limit, the following system holds in $\cD'(\Om)$:
\begin{eqnarray}
  \dx \big( \int_0^{h} \vrh v \, dZ\big) & = & 0, \label{21}\\
  -\dZ \big ( \mu(\vrh) \dZ v\big ) + \dx P(\vrh) & = & 0, \label{22}\\
  \dZ P(\vrh) & = & 0. \label{23}
\end{eqnarray}
Moreover, the horizontal velocity~$v$ satisfies the boundary conditions $v|_{Z=0}=V$, $v|_{Z=h(x)}=0$ and $v$ and $\vrh$ are periodic with respect to the $x$ variable.
\ENDTHEOREM

\REMARK
\begin{enumerate}
\item The boundary conditions on the density are not conserved through the limit $\vep\rightarrow 0$ and this is essential to get a non constant density $\vrh$ in the limit system.
\item The limit model is not coupled with the vertical velocity any more. However, we have a limit equality for the weak limit $w$ of $\vep w_\vep$: $\dZ(\vrh w)=0$.
\item Notice that if $P'$ is not zero  almost everywhere then $\dZ P(\vrh) = 0 \Longleftrightarrow \dZ \vrh = 0$. 
As in the incompressible case, we can integrate twice the first equation~\eqref{22} with respect to variable $Z$ and use the condition~\eqref{21} in order to obtain:
\begin{equation}\label{Reynolds}
\pa_x \left(\frac{h^3}{12} \frac{\vrh P'(\vrh)}{\mu(\vrh)} \pa_x \vrh \right) = \pa_x \left( \frac{\rho \, h}{2} V \right).
\end{equation}
\item It is important to notice that for such a Reynolds equation, maximum principle is proved (see for instance~\cite{Jai}).
Consequently, there exists a constant $\rho_{\text{min}}>0$ such that the solution~$\rho$ of the Reynolds equation~\eqref{Reynolds} satisfies $\rho\geq \rho_{\text{min}}$.
We deduce that for~$\eps$ small enough, we have $\rho_\eps\geq \frac{\rho_{\text{min}}}{2}>0$, that imply that assumptions~\eqref{6} and~\eqref{9} are useless.
\item The proof presented here can easily be extended to the nonstationary case (by using the result of the Theorem~\ref{th-non-stationary}). We thus justify the nonstationary compressible Reynolds equation
\begin{equation*}
\partial_t (\rho \, h) + \pa_x \left(\frac{h^3}{12} \frac{\vrh P'(\vrh)}{\mu(\vrh)} \pa_x \vrh \right) = \pa_x \left( \frac{\rho\, h}{2} V \right),
\end{equation*}
as the limit of the compressible Navier-Stokes equations in thin domain.
\end{enumerate}
\ENDREMARK

\SECTION{5}{Proof of Theorem \statementref{4}}\label{sec-5}  

For all $\eps>0$, the existence of a suitable solution of \eqref{18}--\eqref{20} is given by Theorem \ref{th-stationary}, coupled with classical results in periodical cases (see~\cite{BD}). So, let us consider $(\vrh_\vep,\u_\vep)$ such a solution.\par
The aim is to obtain estimates of $(\vrh_\eps,\u_\eps)$ which are on the one hand non-depending on~$\eps$, on the other hand sufficient to pass to the limit $\eps\rightarrow 0$ in the equations~\eqref{18}--\eqref{20}.
The main difficulty comes from the non-linearities which require strong convergences of some terms.

\SUBSECTION{5.1}{A priori estimates}

To write energy estimates, we take again the method of previous sections. Recall that this method requires the introduction of a velocity lift~$\widetilde{\u}_\vep$.
Within the framework which interests us here (see Figure page~\pageref{fig-domain}), the velocity lift that we will use is:~$\widetilde{\u}_\vep=(\widetilde{v}_\vep,\widetilde{w}_\vep)$ with
\begin{equation*}
  \widetilde{v}_\vep = \left\{ 
  \begin{aligned}
    V(1-Z) \quad &\text{if $0<Z<1$}\\
    0 \hspace{1cm} &\text{if $Z>1$}
  \end{aligned}
\right. \quad \text{and} \quad \widetilde{w}_\vep=0.
\end{equation*}
Note that it is very easy to regularize this velocity field, keeping the form $\widetilde{\u}_\eps=(\widetilde{v}_\eps(Z),0)$.
\LEMMA{5}\label{lemma}
The energy and BD formula write as follows
\begin{equation}
\begin{aligned}\label{thin-energy}
& 2\int_\Om\mu(\vrh_\vep)|\pa_x v_\vep|^2+2\int_\Om\mu(\vrh_\vep)|\pa_Z w_\vep|^2+\int_\Om\mu(\vrh_\vep)\Big |\vep\pa_x w_\vep+\frac{1}{\vep}\pa_Z v_\vep\Big |^2 \\
& \hspace{3cm} + \int_\Om\la(\vrh_\vep)|\pa_x v_\vep+\pa_Z w_\vep|^2
  + r_0\int_\Om \vrh_\vep (v_\vep^2+\vep^2 w_\vep^2)^{3/2}
  = S_1^\vep
\end{aligned}
\end{equation}
\begin{equation}\label{thin-BD}
\begin{aligned}
& \int_\Om\mu(\vrh_\vep)\Big |\vep\pa_x w_\vep-\frac{1}{\vep}\pa_Z v_\vep\Big |^2
+\frac{1}{\vep^2}\frac{C_0}{c_2 M^2}\int_\Om|\pa_x(\xi(\vrh_\vep))^M|^2
+\frac{1}{\vep^4}\frac{C_0}{c_2 M^2}\int_\Om|\pa_Z(\xi(\vrh_\vep))^M|^2 \\
& \hspace{1cm} + r_0\int_\Om \vrh_\vep (v_\vep^2+\vep^2 w_\vep^2)^{3/2}
  + \frac{1}{\vep^2}\frac{C_0 a\ga}{2N^2}\int_\Om |\pa_x (\vrh_\vep^N)|^2
  + \frac{1}{\vep^4}\frac{C_0 a\ga}{2N^2}\int_\Om |\pa_Z (\vrh_\vep^N)|^2 \leq |S_2^\vep|
\end{aligned}
\end{equation}
where 
\begin{equation*}
\begin{aligned}
& S_1^\vep = T_1 + \frac{1}{\vep}\int_\Om \mu(\vrh_\vep)\, \Big(\frac{1}{\vep} \,\pa_Z v_\vep + \vep \, \pa_x w_\vep \Big )\, \dz \widetilde{v_\eps}
 + r_0 \int_\Om \vrh_\vep\, \widetilde{v}_\vep \,  v_\eps \, (v_\vep^2+\vep^2 w_\vep^2)^{1/2}, \\
& \text{with} \quad T_1 = \int_\Om\vrh_\vep v_\vep\pa_x v_\vep \, \widetilde{v}_\vep + \int_\Om\vrh_\vep w_\vep\pa_Z v_\vep \, \widetilde{v}_\vep,\\
& S_2^\vep = S_1^\vep + 2 \int_\Omega \dx \mu(\vrh_\vep) \, w_\eps \, \dZ \widetilde v_\eps + 2 \, r_0 \int_\Omega \big (v_\eps \, \dx \mu(\vrh_\vep) + w_\eps \, \dZ \mu(\vrh_\vep)\big )\, (v_\eps^2+\eps^2 \, w_\eps^2)^{1/2}.
\end{aligned}
\end{equation*}

\ENDLEMMA

{\bf Proof of the energy estimate~\eqref{thin-energy}}\par
We first note that the function $\widetilde{\u}_\vep$ permits to get homogeneous Dirichlet boundary conditions on the corrected velocity $\u_\vep-\widetilde{\u}_\vep$.
To get a first energy identity, we sum equation \eqref{19} multiplied by the new horizontal velocity $v_\vep-\widetilde{v}_\vep$ and equation \eqref{20} multiplied by $\vep(w_\vep-\widetilde w_\vep)=\vep w_\vep$.
We obtain
\begin{equation}\label{estimate1045}
\begin{aligned}
& \int_\Om\vrh_\vep v_\vep\pa_x v_\vep (v_\vep-\widetilde{v}_\vep)+ \int_\Om\vrh_\vep w_\vep\pa_Z v_\vep (v_\vep-\widetilde{v}_\vep)
+\vep^2\int_\Om\vrh_\vep v_\vep w_\vep\pa_x w_\vep +\vep^2\int_\Om\vrh_\vep w_\vep^2\pa_Z w_\vep \\
& \quad = 2\int_\Om\pa_x(\mu(\vrh_\vep) \pa_x v_\vep)(v_\vep-\widetilde{v}_\vep)
+ \frac{1}{\vep^2}\int_\Om\pa_Z\big (\mu(\vrh_\vep) \pa_Z v_\vep\big )(v_\vep-\widetilde{v}_\vep)+ \int_\Om\pa_Z\big (\mu(\vrh_\vep) \pa_x w_\vep\big )(v_\vep-\widetilde{v}_\vep) \\
& \hspace{2cm} + \int_\Om\pa_x\big (\mu(\vrh_\vep)\pa_Z v_\vep\big ) w_\vep+ \vep^2\int_\Om\pa_x\big (\mu(\vrh_\vep) \pa_x w_\vep\big ) w_\vep 
+ 2\int_\Om\pa_Z\big (\mu(\vrh_\vep) \pa_Z w_\vep\big ) w_\vep \\
& \hspace{2cm} + \int_\Om\pa_x\big (\lambda(\vrh_\vep) (\pa_x v_\vep + \pa_Z w_\vep)\big )(v_\vep-\widetilde{v}_\vep)+\int_\Om\pa_Z\big (\lambda(\vrh_\vep) (\pa_x v_\vep + \pa_Z w_\vep)\big )w_\vep \\ 
& \hspace{2cm} - \frac{1}{\vep^2} \int_\Om\pa_x P(\vrh_\vep)(v_\vep-\widetilde{v}_\vep)
-\frac{1}{\vep^2} \int_\Om\pa_Z P(\vrh_\vep)w_\vep \\
& \hspace{2cm} - r_0\int_\Om\vrh_\vep |v_\vep^2+\vep^2 w_\vep^2|^{\frac{1}{2}}v_\vep(v_\vep-\widetilde{v}_\vep)-r_0\vep^2\int_\Om\vrh_\vep|v_\vep^2+\vep^2 w_\vep^2|^{\frac{1}{2}} w_\vep^2.
\end{aligned}
\end{equation}
By the complexity of this equation, we will deal with these terms by group.
The terms of the first line (left hand side of the equation) are known as convection terms.
Those of the three following lines will be called viscous terms.
The fifth line corresponds to the terms of pressure whereas the last line contains the friction ones.
\par
{\bf Convection terms -}
We first use integrations by parts (without boundary terms thanks to the $x$-periodicity and the boundary conditions for~$w_\eps$) and the divergence free conditions:
\begin{equation*}
\int_\Om\vrh_\vep v_\vep\pa_x v_\vep \, v_\vep + \int_\Om\vrh_\vep w_\vep\pa_Z v_\vep \, v_\vep 
= - \int_\Om\vrh_\vep v_\vep\, v_\vep \dx v_\vep - \int_\Om\vrh_\vep w_\vep\, v_\vep \dZ v_\vep. 
\end{equation*}
This contribution is zero (since it equals to its opposite). In the same way we have
\begin{equation*}
\vep^2\int_\Om\vrh_\vep v_\vep w_\vep\pa_x w_\vep +\vep^2\int_\Om\vrh_\vep w_\vep^2\pa_Z w_\vep = 0.
\end{equation*}
Finally, all the terms of the left hand side of~\eqref{estimate1045} disappear except those containing~$\widetilde{v_\eps}$, denoted~$T_1$:
\begin{equation*}
T_1 = \int_\Om\vrh_\vep v_\vep\pa_x v_\vep \, \widetilde{v}_\vep + \int_\Om\vrh_\vep w_\vep\pa_Z v_\vep \, \widetilde{v}_\vep.
\end{equation*}
{\bf Viscous terms -}
These terms are easily computed using integrations by parts.
In any integration by parts, no boundary integral term appear thanks to the vertical correction induced by $\widetilde{\u}_\vep$ and the periodicity in~$x$.
\par
{\bf Pressure terms -}
We also remark that, using $\div(\vrh_\vep \u_\vep)=\pa_x(\vrh_\vep v_\vep)+\pa_Z(\vrh_\vep w_\vep)=0$ and $\dx \widetilde{v_\vep}=0$, the pressure contributions vanish.
In fact, noting $\Pi'(\vrh_\vep)=\frac{P'(\vrh_\vep)}{\vrh_\vep}$, we have
\[
\int_\Om\pa_x P(\vrh_\vep)(v_\vep-\widetilde{v}_\vep)
+ \int_\Om\pa_Z P(\vrh_\vep)w_\vep 
= -\int_\Om \Pi(\vrh_\vep) \big (\pa_x (\vrh_\vep v_\vep)+\pa_Z(\vrh_\vep w_\vep)\big )-\int_\Om P(\vrh_\vep)\pa_x \widetilde{v}_\vep
= 0.
\]
{\bf Friction terms -}
Clearly , the friction terms (that is the terms containing the friction coefficient~$r_0$) appearing in estimate~\eqref{estimate1045} write
\begin{equation*}
-r_0\int_\Om\vrh_\vep |v_\vep^2+\vep^2 w_\vep^2|^{\frac{3}{2}} + r_0\int_\Om\vrh_\vep |v_\vep^2+\vep^2 w_\vep^2|^{\frac{1}{2}}v_\vep\,\widetilde{v}_\vep.
\end{equation*}
All these calculations give the first identity of Lemma~\statementref{5}.
\par\vspace{0.3cm}
{\bf Proof of the energy estimate~\eqref{thin-BD}}\par
Refering to Section \sectionref{2}, we can also write another equality related to the BD entropy.
To get this BD formula, we sum equation \eqref{BD-eps1} multiplied by $v_\vep-\widetilde{v}_\vep+V_\vep$ and equation \eqref{BD-eps2} multiplied by $w_\vep-\widetilde{w}_\vep+W_\vep$ (for information, recall that $V_\vep = \frac{2\pa_x\mu(\vrh_\vep)}{\vrh_\vep}$ and $W_\vep=\frac{2\pa_Z\mu(\vrh_\vep)}{\vep\vrh_\vep}$). We obtain
\begin{equation}\label{estimate1149}
\begin{aligned}
& \int_\Om\vrh_\vep v_\vep\pa_x(v_\vep+V_\vep)(v_\vep-\widetilde{v}_\vep+V_\vep) + \int_\Om\vrh_\vep w_\vep\pa_Z(v_\vep+V_\vep)(v_\vep-\widetilde{v}_\vep+V_\vep) \\
& + \int_\Om\vrh_\vep v_\vep \pa_x(\vep w_\vep+W_\vep)(\vep w_\vep+W_\vep) + \int_\Om\vrh_\vep w_\vep \pa_Z(\vep w_\vep+W_\vep)(\vep w_\vep+W_\vep)\\
& \qquad = \int_\Om\pa_Z(\mu(\vrh_\vep) \pa_x w_\vep)(v_\vep-\widetilde{v}_\vep+V_\vep) - \frac{1}{\vep^2}\int_\Om\pa_Z(\mu(\vrh_\vep) \pa_Z v_\vep)(v_\vep-\widetilde{v}_\vep+V_\vep) \\
& \hspace{2cm} - \vep\int_\Om\pa_x(\mu(\vrh_\vep) \pa_x w_\vep)(\vep w_\vep+W_\vep) + \frac{1}{\vep}\int_\Om\pa_x(\mu(\vrh_\vep) \pa_Z v_\vep) (\vep w_\vep+W_\vep) \\
& \hspace{2cm} - \frac{1}{\vep^2}\int_\Om\pa_x P(\vrh_\vep)(v_\vep-\widetilde{v}_\vep+V_\vep) - \frac{1}{\vep^3}\int_\Om\pa_Z P(\vrh_\vep)(\vep w_\vep+W_\vep)\\
& \hspace{2cm} - r_0\int_\Om\vrh_\vep|v_\vep^2+\vep^2 w_\vep^2|^{\frac{1}{2}} v_\vep(v_\vep-\widetilde{v}_\vep+V_\vep) - r_0\int_\Om\vep\vrh_\vep|v_\vep^2+\vep^2 w_\vep^2|^{\frac{1}{2}} w_\vep (\vep w_\vep+W_\vep).
\end{aligned}
\end{equation}
In addition to those which are common with the energy, we deal with every terms:
\par
{\bf Convection terms -}
Let us look at the first line of the equation~\eqref{estimate1149}.
The terms which are not already treated to obtain identity~\eqref{thin-energy} are the following ones:
\begin{equation*}
\begin{aligned}
& \int_\Omega \rho_\eps \, v_\eps \, \dx v_\eps \, V_\eps
+ \int_\Omega \rho_\eps \, v_\eps \, \dx V_\eps \, (v_\eps-\widetilde v_\eps)
+ \int_\Omega \rho_\eps \, v_\eps \, \dx V_\eps \, V_\eps \\
& \qquad + \int_\Omega \rho_\eps \, w_\eps \, \dZ v_\eps \, V_\eps
+ \int_\Omega \rho_\eps \, w_\eps \, \dZ V_\eps \, (v_\eps-\widetilde v_\eps)
+ \int_\Omega \rho_\eps \, w_\eps \, \dZ V_\eps \, V_\eps.
\end{aligned}
\end{equation*}
The divergence free condition $\pa_x(\vrh_\vep v_\vep)+\pa_Z(\vrh_\vep w_\vep)=0$ is also strongly used here.
With the boundary conditions (periodicity with respect to~$x$ and boundary conditions on~$w$ for $Z\in \{0,1\}$) we deduce that all terms are equal to zero by integration by parts, except this one: $\dsp \int_\Omega \rho_\eps \, w_\eps \, V_\eps \, \dZ \widetilde v_\eps$.
\par
Since $\rho_\eps V_\eps = 2 \dx \mu(\rho_\eps)$, we write this additive term as $\dsp 2 \int_\Omega \dx \mu(\vrh_\vep) \, w_\eps \, \dZ \widetilde v_\eps$.
\par
{\bf Viscous terms -}
The only additive terms compared the first energy identity (see equation~\eqref{thin-energy}) are the following
\begin{equation*}
\begin{aligned}
  \int_\Om\pa_Z(\mu(\vrh_\vep) \pa_x w_\vep)V_\vep 
& - \frac{1}{\vep^2}\int_\Om\pa_Z(\mu(\vrh_\vep) \pa_Z v_\vep)V_\vep
  - \vep\int_\Om\pa_x(\mu(\vrh_\vep) \pa_x w_\vep)W_\vep + \frac{1}{\vep}\int_\Om\pa_x(\mu(\vrh_\vep) \pa_Z v_\vep)W_\vep \\
& = \int_\Om\mu(\vrh_\vep) \pa_x w_\vep(\underbrace{\vep\pa_x W_\vep-\pa_Z V_\vep}_{=0})
  + \frac{1}{\vep^2}\int_\Om\mu(\vrh_\vep) \pa_Z v_\vep(\underbrace{\pa_Z V_\vep-\vep\pa_x W_\vep}_{=0}) = 0.
\end{aligned}
\end{equation*}
The equality $\vep\pa_x W_\vep-\pa_Z V_\vep=0$ comes from to the following computation
\begin{equation*}
\begin{aligned}
\vep\pa_x W_\vep-\pa_Z V_\vep 
& = 2 \dx \Big( \frac{\dZ \mu(\vrh_\vep)}{\rho_\vep} \Big) - 2 \dZ \Big( \frac{\dx \mu(\vrh_\vep)}{\rho_\vep} \Big)
= 2 \dZ \mu(\vrh_\vep) \, \dx \Big( \frac{1}{\rho_\vep} \Big) - 2 \dx \mu(\vrh_\vep) \, \dZ \Big( \frac{1}{\rho_\vep} \Big) \\
& = 2 \mu'(\vrh_\vep) \Big[ \dZ \rho_\vep \, \dx \Big( \frac{1}{\rho_\vep} \Big) - \dx \rho_\vep \, \dZ \Big( \frac{1}{\rho_\vep} \Big) \Big]
= -\frac{2 \mu'(\vrh_\vep)}{\rho_\vep^2} \Big[ \dZ \rho_\vep \, \dx \rho_\vep - \dx \rho_\vep \, \dZ \rho_\vep \Big]
= 0.
\end{aligned}
\end{equation*}
\par
{\bf Pressure terms -}
For the pressure terms, we can rewrite what we wrote for the energy, say, the tests against $\u_\vep-\widetilde{\u}_\vep$ are equal to zero thanks to \eqref{18} and $\pa_x v_\vep=0$.
The remaining part of the pressure contributions have been discussed and expressed in~\eqref{14}, let's recall it here:
\begin{equation*}
\begin{aligned}
& \frac{1}{\vep^2}\int_\Om\pa_x P(\vrh_\vep)V_\vep
+ \frac{1}{\vep^3}\int_\Om\pa_Z P(\vrh_\vep)W_\vep \\
& \qquad = \frac{2}{\vep^2}\int_\Om\pa_x \Big (p_h(\vrh_\vep)+p_c(\vrh_\vep)\Big )\frac{\pa_x\mu(\vrh_\vep)}{\vrh_\vep}
+ \frac{2}{\vep^4}\int_\Om\pa_Z \Big (p_h(\vrh_\vep)+p_c(\vrh_\vep)\Big )\frac{\pa_Z\mu(\vrh_\vep)}{\vrh_\vep} \\
&  \qquad \qquad \geq 
  \frac{1}{\vep^2}\frac{c_0'}{c_2 M^2}\int_\Om|\pa_x(\xi(\vrh_\vep))^M|^2
+ \frac{1}{\vep^4}\frac{c_0'}{c_2 M^2}\int_\Om|\pa_Z(\xi(\vrh_\vep))^M|^2 \\
& \qquad \qquad  \qquad \qquad + \frac{1}{\vep^2}\frac{C_0 a\ga}{2N^2}\int_\Om |\pa_x (\vrh_\vep^N)|^2
+ \frac{1}{\vep^4}\frac{C_0 a\ga}{2N^2}\int_\Om |\pa_Z (\vrh_\vep^N)|^2.
\end{aligned}
\end{equation*}
\par
{\bf Friction terms -}
Let's now deal with the friction terms, in addition to those which are common with the energy, we have to say some words about the ones which are specific for the BD formula.
Using the definition of~$V_\eps$ and~$W_\eps$, these additive terms are written
\begin{equation*}
\begin{aligned}
r_0 \int_\Omega \rho_\eps \, (v_\eps^2+\eps^2 \, w_\eps^2)^{1/2} \, (v_\eps \, V_\eps + \eps\, w_\eps \, W_\eps)
= 2 \, r_0 \int_\Omega (v_\eps^2+\eps^2 \, w_\eps^2)^{1/2} \, \big (v_\eps \, \dx \mu(\vrh_\vep) + w_\eps \, \dZ \mu(\vrh_\vep)\big ).
\end{aligned}
\end{equation*}
\par
Putting all these inequalities together, we obtain~\eqref{thin-BD}.
\endsquare
\par\vspace{0.5cm}
The left hand sides of the estimates~\eqref{thin-energy} and~\eqref{thin-BD} make appear the $L^1(\Omega)$-norm of the following terms
\begin{equation*}
\begin{aligned}
& \mu_\vep | \dx v_\eps |^2, \quad
\mu_\vep | \dZ w_\eps |^2, \quad
\eps^2 \mu_\vep | \dx w_\eps |^2, \quad
\frac{\mu_\vep}{\eps^2} | \dZ v_\eps |^2, \quad
\lambda_\vep | \dx v_\eps |^2, \quad
\lambda_\vep | \dZ w_\eps |^2, \quad
r_0 \rho_\vep | v_\eps |^3, \quad
\eps^3 r_0 \rho_\vep | w_\eps |^3, \\
& \eps^2 \mu_\vep | \dx w_\eps |^2, \quad
\frac{\mu_\vep}{\eps^2} | \dZ v_\eps |^2, \quad
\frac{1}{\eps^2} | \dx (\xi_\vep^M) |^{2}, \quad
\frac{1}{\eps^4} | \dZ (\xi_\vep^M) |^{2}, \quad
\frac{1}{\eps^2} | \dx (\rho_\vep^N) |^2 \quad \text{and} \quad
\frac{1}{\eps^4} | \dZ (\rho_\vep^N) |^2,
\end{aligned}
\end{equation*}
where we have noted $\mu_\vep=\mu(\vrh_\vep)$, $\la_\vep=\la(\vrh_\vep)$ and $\xi_\vep=\xi(\vrh_\vep)$.
We now prove that the right hand sides of the estimates~\eqref{thin-energy} and~\eqref{thin-BD}, that is the quantities~$S_1^\vep$ and~$S_2^\vep$, can be controlled by such terms.
\vskip 10pt
{\bf Control of $S_1^\vep$ -}
>From the definition of~$S_1^\eps$, we express~$S_1^\vep$ as follows: $S_1^\eps = T_1+T_2+T_3$.
\par
The contribution~$T_1$ is written\footnote{Note that this term is not treated as the corresponding one in the proof of Theorem~\ref{th-stationary}.
In the proof of Theorem~\ref{th-stationary}, we write $T_1=-\int_\Omega \rho_\eps \, v_\eps \, w_\eps \ \dZ \widetilde{v_\eps}$. This is correct but inappropriate in our situation.
In fact, we want to obtain estimates with respect to the parameter~$\eps$ and we have no control on $w_\eps$, but only on $\eps\, w_\eps$ and on $\frac{1}{\eps}\dZ v_\eps$.}
\begin{equation*}
\begin{aligned}
T_1 
& = \int_\Omega \rho_\eps \, v_\eps \, \dx v_\eps \, \widetilde{v_\eps}
+ \int_\Omega \rho_\eps \, w_\eps \, \dZ v_\eps \, \widetilde{v_\eps} \\
& = \int_\Omega \big( \rho_\eps^{1/3} \, v_\eps \big) \, \big( \rho_\eps^{m/2} \, \dx v_\eps \big) \, \big( \rho_\eps^{(4-3m)/6} \widetilde{v_\eps} \big)
+ \int_\Omega \big( \rho_\eps^{1/3} \, \eps\, w_\eps \big) \, \big( \frac{\rho_\eps^{m/2} \, \dZ v_\eps }{\eps}\big) \, \big( \rho_\eps^{(4-3m)/6} \widetilde{v_\eps} \big) \\
& \leq \frac{r_0}{4} \int_\Omega \rho_\eps v_\eps^3 
+ \delta \int_\Omega \rho_\eps^m |\dx v_\eps|^2 
+ \frac{C}{r_0^2\, \delta^3} |\widetilde{v_\eps}|_\infty^6\int_\Omega \rho_\eps^{4-3m} \\
& \hspace{2cm}
+ \frac{r_0\, \eps^3}{4} \int_\Omega \rho_\eps \, w_\eps^3 
+ \frac{\delta}{\eps^2} \int_\Omega \rho_\eps^m |\dZ v_\eps|^2 
+ \frac{C}{r_0^2\, \delta^3} |\widetilde{v_\eps}|_\infty^6\int_\Omega \rho_\eps^{4-3m},
\end{aligned}
\end{equation*}
In this last inequality,~$\delta$ can be choosen as small as possible (due to the Young inequality).
Moreover the constant~$C$ does not depend on the physical constants~$\eps$, $r_0$, $c_1$... nor~$\delta$.
The term $\delta \int_\Omega \rho_\eps^m |\dx v_\eps|^2$ is bounded using the assumptions~\eqref{7} on $\mu$ as follows:
\begin{equation*}
\begin{aligned}
\delta \int_\Omega \rho_\eps^m |\dx v_\eps|^2  
& = \delta \int_\Omega \rho_\eps^m |\dx v_\eps|^2 \, \mathds{1}_{\rho>A}
+ \delta \int_\Omega \rho_\eps^m |\dx v_\eps|^2\, \mathds{1}_{\rho<A} \\
& \leq \frac{\delta}{c_1} \int_\Omega \mu(\vrh_\vep) \, |\dx v_\eps|^2 \, \mathds{1}_{\rho>A}
+ \frac{\delta}{A^{m-n}} \int_\Omega \rho_\eps^n |\dx v_\eps|^2\, \mathds{1}_{\rho<A} \\
& \leq \frac{\delta}{c_1} \int_\Omega \mu(\vrh_\vep) \, |\dx v_\eps|^2 \, \mathds{1}_{\rho>A}
+ \frac{\delta}{c_0\, A^{m-n}} \int_\Omega \mu(\vrh_\vep) \, |\dx v_\eps|^2\, \mathds{1}_{\rho<A} \\
\end{aligned}
\end{equation*}
Taking $\delta=\min \{ \frac{c_1}{2},\frac{c_0\, A^{m-n}}{2}\}$ we obtain
\begin{equation*}
\delta \int_\Omega \rho_\eps^m |\dx v_\eps|^2  
\leq \frac{1}{2} \int_\Omega \mu(\vrh_\vep) \, |\dx v_\eps|^2.
\end{equation*}
In the same way, we estime $\frac{\delta}{\eps^2} \int_\Omega \rho_\eps^m |\dZ v_\eps|^2$ and we obtain (note also that $|\widetilde{v_\eps}|_\infty$ is bounded by~$1$)
\begin{equation*}
T_1 
\leq \frac{r_0}{4} \int_\Omega \rho_\eps v_\eps^3 
+ \frac{1}{2} \int_\Omega \mu(\vrh_\vep) \, |\dx v_\eps|^2
+ \frac{1}{2\, \eps^2} \int_\Omega \mu(\vrh_\vep) \, |\dZ v_\eps|^2
+ \frac{r_0\, \eps^3}{4} \int_\Omega \rho_\eps \, w_\eps^3
+ \frac{C}{r_0^2\, \delta^3} \int_\Omega \rho_\eps^{4-3m}.
\end{equation*}
Since $|\dZ \widetilde{v_\eps}|\leq V$, the contribution~$T_2$ is controlled as follows
\begin{equation*}
\begin{aligned}
T_2 
& = \frac{1}{\vep}\int_\Om \mu(\vrh_\vep)\, \big(\frac{1}{\vep} \,\pa_Z v_\vep + \vep \, \pa_x w_\vep \big )\, \dz \widetilde{v_\eps}\\
& \leq \frac{V}{\vep} \int_\Om\mu(\vrh_\vep)\big| \frac{1}{\vep} \,\pa_Z v_\vep + \vep \, \pa_x w_\vep \big| \\
& \leq \frac{1}{2}\int_\Om\mu(\vrh_\vep)\big| \frac{1}{\vep} \,\pa_Z v_\vep + \vep \, \pa_x w_\vep \big|^2+\frac{V^2}{2\vep^2}\int_\Om \mu(\vrh_\vep).
\end{aligned}
\end{equation*}
Using assumptions~\eqref{6} and~\eqref{7}, we deduce that
\begin{equation*}
T_2 \leq \frac{1}{2}\int_\Om\mu(\vrh_\vep)\big| \frac{1}{\vep} \,\pa_Z v_\vep + \vep \, \pa_x w_\vep \big|^2 + \frac{V^2}{2c_0\, \vep^2} + \frac{V^2}{2c_1\, \vep^2}\int_\Om \vrh_\vep^m.
\end{equation*}
The term~$T_3$ is treated using the Young inegality (in the following form: $A B^2 \leq \delta A B^3 + \frac{4}{27 \delta^2}A$, for all $\delta>0$). We obtain
\begin{equation*}
  \begin{aligned}
T_3 
& = r_0 \int_\Om \vrh_\vep\, \widetilde{v}_\vep \,  v_\eps \, (v_\vep^2+\vep^2 w_\vep^2)^{1/2} \\
& \leq r_0 \int_\Om \vrh_\vep \,  |v_\eps |^2 + r_0 \int_\Om \vrh_\vep\,  |v_\eps| \, |\vep \, w_\vep| \\
& \leq \frac{3r_0}{2} \int_\Om \vrh_\vep \,  |v_\eps |^2 + \frac{r_0}{2} \int_\Om \vrh_\vep \,  |\eps\, w_\eps |^2 \\
& \leq \frac{r_0}{4} \int_\Omega \rho_\eps \, |v_\eps|^3 
+ \frac{r_0\, \eps^3}{4} \int_\Omega \rho_\eps \, |w_\eps|^3
+ C \, r_0 \int_\Omega \rho_\eps. 
  \end{aligned}
\end{equation*}
Here, the constant~$C$ does not depend on~$\eps$ nor on~$r_0$.
\par
{\bf Control of $S_2^\vep$ -}
>From the definition of~$S_2^\eps$, we write~$S_2^\vep$ as follows: $S_2^\eps = S_1^\eps + T_4+T_5$.
\par
Using an integration by part and $|\dZ \widetilde{v_\eps}|\leq V$, the contribution~$T_4$ is controlled by
\[
T_4 
\leq 2 V \int_\Om\mu(\vrh_\vep) |\pa_x w_\vep|
\leq \frac{2V^2}{\vep^2}\int\mu(\vrh_\vep)
+ \frac{\vep^2}{2}\int_\Om\mu(\vrh_\vep)|\pa_x w_\vep|^2.
\]
Using assumptions~\eqref{6} and~\eqref{7}, we deduce that
\begin{equation*}
T_4 
\leq \frac{\vep^2}{2}\int_\Om\mu(\vrh_\vep)|\pa_x w_\vep|^2
+ \frac{2V^2}{c_0\, \vep^2} + \frac{2V^2}{c_1\, \vep^2}\int_\Om \vrh_\vep^m.
\end{equation*}

Finally, the term~$T_5$ is written
\begin{equation*}
T_5 = 2 \, r_0 \int_\Omega (v_\eps \, \dx \mu(\vrh_\vep) + w_\eps \, \dZ \mu(\vrh_\vep))\, (v_\eps^2+\eps^2 \, w_\eps^2)^{1/2}.
\end{equation*}
For sake of simplicity, we only treat one example of this contribution (more precisely, the term $2\, r_0 \int_\Omega v_\eps^2 \, \dx \mu(\vrh_\vep) $), the other terms are similar.
We have
\begin{equation*}
\begin{aligned}
2\, r_0 \int_\Omega v_\eps^2 \, \dx \mu(\vrh_\vep) 
& = -4\, r_0 \int_\Omega \mu(\vrh_\vep)\, v_\eps\, \dx v_\eps
\leq \frac{1}{2} \int_\Omega \mu(\vrh_\vep) \, |\dx v_\eps|^2 + 8\,r_0^2 \int_\Omega \mu(\vrh_\vep) \, |v_\eps|^2 \\
& \leq \frac{1}{2} \int_\Omega \mu(\vrh_\vep) \, |\dx v_\eps|^2 
+ 8\,r_0^2 \int_\Omega \mu(\vrh_\vep) \, |v_\eps|^2 \, \mathds{1}_{\rho<A}
+ 8\,r_0^2 \int_\Omega \mu(\vrh_\vep) \, |v_\eps|^2 \, \mathds{1}_{\rho>A}.
\end{aligned}
\end{equation*}
Using succesively the assumption~\eqref{6}, the fact that $n>2/3$ and the Young inequality we obtain
\begin{equation*}
8\,r_0^2 \int_\Omega \mu(\vrh_\vep) \, |v_\eps|^2 \, \mathds{1}_{\rho<A} 
\leq \frac{8\, r_0^2}{c_0} \int_\Omega \rho^n \, |v_\eps|^2\, \mathds{1}_{\rho<A} 
\leq \frac{8\, r_0^2\, A^{n-\frac{2}{3}}}{c_0} \int_\Omega \rho^{2/3} \, |v_\eps|^2
\leq \frac{r_0}{4} \int_\Omega \rho_\eps \, |v_\eps|^3 + \frac{C\, r_0^4\, A^{3n-2}}{c_0^3}.
\end{equation*}
In the same way, we use the assumption~\eqref{7} and the Young inequality to obtain
\begin{equation*}
8\,r_0^2 \int_\Omega \mu(\vrh_\vep) \, |v_\eps|^2 \, \mathds{1}_{\rho>A}
\leq \frac{8\, r_0^2}{c_1} \int_\Omega \rho^m \, |v_\eps|^2
\leq \frac{r_0}{4} \int_\Omega \rho_\eps \, |v_\eps|^3 + \frac{C\, r_0^4}{c_1^3} \int_\Omega \rho_\eps^{3m-2}.
\end{equation*}
Here again, the constant~$C$ does not depend on~$\eps$, $r_0$, $c_0$ or~$c_1$.
\par
Finally, the term $T_5$ is controlled by terms like
\begin{equation*}
T_5 
\leq \frac{1}{2} \int_\Omega \mu(\vrh_\vep) \, |\dx v_\eps|^2 
+ \frac{r_0}{4} \int_\Omega \rho_\eps \, |v_\eps|^3
+ \frac{C\, r_0^4}{c_1^3} \int_\Omega \rho_\eps^{3m-2}
+ \frac{C\, r_0^4\, A^{3n-2}}{c_0^3}.
\end{equation*}
\vskip 5pt
{\bf Estimates -}
The result of Lemma~\ref{lemma} and the control of the terms~$S_1^\vep$ and~$S_2^\vep$ allow to obtain the following estimates (note that in this estimate and in the following ones, the only constant that we will reveal will be~$\eps$, the others will be taken equal to~$1$ to simplify calculations).
\begin{equation}
\begin{aligned}\label{thin-estimate}
& \int_\Om \mu(\vrh_\vep) \, |\pa_x v_\vep|^2
+ \int_\Om \mu(\vrh_\vep) \, |\pa_Z w_\vep|^2
+ \eps^2 \int_\Om \mu(\vrh_\vep) \, |\dx w_\vep|^2
+ \frac{1}{\eps^2} \int_\Om \mu(\vrh_\vep) \, |\dZ v_\vep|^2
+ \int_\Omega \rho_\eps \, |v_\eps|^3 \\
&
+ \eps^3 \int_\Omega \rho_\eps \, |w_\eps|^3
+ \frac{1}{\vep^2} \int_\Om |\pa_x(\xi_\eps^M)|^2
+ \frac{1}{\vep^4} \int_\Om |\pa_Z(\xi_\eps^M)|^2
+ \frac{1}{\vep^2} \int_\Om |\pa_x(\vrh_\vep^N)|^2
+ \frac{1}{\vep^4} \int_\Om |\pa_Z(\vrh_\vep^N)|^2 \\
& \hspace{2cm}
\leq
\int_\Omega \rho_\eps^{4-3m}
+ \frac{1}{\eps^2}
+ \frac{1}{\eps^2} \int_\Omega \rho_\eps^m
+ \int_\Omega \rho_\eps
+ \int_\Omega \rho_\eps^{3m-2}
+ 1.
\end{aligned}
\end{equation}
\par
The right hand side of the estimate~\eqref{thin-estimate} make appear different powers of $\rho_\eps$.
We must show that all these terms can be absorbed by some left hand side terms, taking~$\eps$ small enough.
\par
Since $m>1$, we have $3m-2>m>1$ and $3m-2>4-3m$.
So, we must control the term $\int_\Omega \rho_\eps^{3m-2}$ and the term $\int_\Omega \rho_\eps^{4-3m}$ when $4-3m<0$.
\par\vspace{0.3cm}
$\bullet$
With the Poincaré inequality, a control on~$|\nabla (\rho_\eps^{N})|_{L^2(\Omega)}$ implies a control on~$|\rho_\eps^{N}|_{H^1(\Omega)}$.
Due to the Sobolev embeddings, this allows to control~$|\rho_\eps^{N}|_{L^q(\Omega)}$ for all $q\geq \frac{2d}{d-2}$.
\par
Consequently, if 
\begin{equation}\label{C1tilde}\tag{$\widetilde{C1}$}
3m-2 \leq qN
\end{equation}
then the term $\int_\Omega \rho_\eps^{3m-2}$ can be absorbed by the terms $\frac{1}{\vep^4} \int_\Omega |\dZ (\vrh_\vep^N)|^2$ and $\frac{1}{\vep^2} \int_\Omega |\dx (\vrh_\vep^N)|^2$ as soon as~$\eps$ is small enough\footnote{If we assume that $3m-2 < qN$, that is the condition~\eqref{C1}, then it is possible to control $\int_\Omega \rho_\eps^{3m-2}$ without taking $\eps$ small but just using a Young inequality, see for instance the proof of Theorem~\ref{th-stationary} where such a method is used.}.
\par\vspace{0.3cm}
$\bullet$ Using the same arguments, from the definition of the function~$\xi$, the term $\int_\Omega |\pa_Z(\xi_\vep^M)|^2$ allows to control $\int_\Omega \rho_\eps^{qM}$, and recalling that $M<0$ we deduce a control of $1/\rho_\eps$ in $L^{-qM}(\Omega)$.
Consequently, when $4-3m<0$ we can absorbed the term $\int_\Omega \rho_\eps^{4-3m}$ by the term $\frac{1}{\vep^4}\int_\Om|\pa_Z(\xi_\vep^M)|^2 $ as soon as~$\eps$ is small enough, under the condition
\begin{equation}\label{C3}\tag{$C3$}
qM \leq 4-3m.
\end{equation}
\par\vspace{0.3cm}
Finally, under the assumptions~$3m-2 \leq qN$ and~$qM \leq 4-3m$ we obtain for~$\eps$ small enough:
\begin{equation*}
\begin{aligned}
& \int_\Om \mu(\vrh_\vep) \, |\pa_x v_\vep|^2
+ \int_\Om \mu(\vrh_\vep) \, |\pa_Z w_\vep|^2
+ \eps^2 \int_\Om \mu(\vrh_\vep) \, |\dx w_\vep|^2
+ \frac{1}{\eps^2} \int_\Om \mu(\vrh_\vep) \, |\dZ v_\vep|^2
+  \int_\Omega \rho_\eps \, |v_\eps|^3 \\
& 
+ \eps^3 \int_\Omega \rho_\eps \, |w_\eps|^3
+ \frac{1}{\vep^2} \int_\Om|\pa_x (\xi(\vrh_\vep)^M)|^2
+ \frac{1}{\vep^4} \int_\Om|\pa_Z (\xi(\vrh_\vep)^M)|^2
+ \frac{1}{\vep^2} \int_\Om |\pa_x (\vrh_\vep^{N})|^2
+ \frac{1}{\vep^4} \int_\Om |\pa_Z (\vrh_\vep^{N})|^2 
\leq \frac{1}{\eps^2}.
\end{aligned}
\end{equation*}
We deduce the following bounds:
\begin{eqnarray}
  \|\sqrt{\mu(\vrh_\vep)}\pa_x v_\vep\|_{L^2(\Om)} & \leq & \nfrac{1}{\vep} \label{25} \\
  \|\sqrt{\mu(\vrh_\vep)}\pa_Z v_\vep\|_{L^2(\Om)} & \leq & 1 \label{26} \\
  \|\sqrt{\mu(\vrh_\vep)}\pa_x w_\vep\|_{L^2(\Om)} & \leq & \nfrac{1}{\vep^2} \label{27} \\
  \|\sqrt{\mu(\vrh_\vep)}\pa_Z w_\vep\|_{L^2(\Om)} & \leq & \nfrac{1}{\vep} \label{28} \\
  \|\pa_x (\xi(\vrh_\vep)^M) \|_{L^2(\Om)} & \leq & 1 \label{29} \\
  \|\pa_Z (\xi(\vrh_\vep)^M) \|_{L^2(\Om)} & \leq & \vep \label{30} \\
  \|\pa_x (\vrh_\vep^{N})\|_{L^2(\Om)} & \leq & 1 \label{31} \\
  \|\pa_Z (\vrh_\vep^{N})\|_{L^2(\Om)} & \leq & \vep \label{32} \\
  \|\sqrt{\vrh_\vep}v_\vep^{\frac{3}{2}}\|_{L^2(\Om)} & \leq & \nfrac{1}{\vep}\label{33a} \\
  \|\sqrt{\vrh_\vep}w_\vep^{\frac{3}{2}}\|_{L^2(\Om)} & \leq & \nfrac{1}{\vep^{\frac{5}{2}}}\label{33b}
\end{eqnarray}
where $M=\frac{n-\al-1}{2}<0$ and $N=\frac{n+\gamma-1}{2}>0$.

\SUBSECTION{5.2}{Compactness on the density}

Let us recall that, for simplicity of the notations, the preceding calculations were carried out on a two dimensional domain, i.e. on $\Omega\subset \R^d$ with $d=2$.
Of course, the latter remains valid in higher dimension, in particular in dimension $d=3$.
In this paragraph, we will strongly use Sobolev injections which are depending on the dimension.
In order to cover the general case $d\in\{2,3\}$, we will note~$q$ any real such that $H^1(\Om)\subset L^q(\Om)$ with continuous injection.
\par
The previous estimates show that $\vrh_\vep^M$ and $\vrh_\vep^N$ are bounded in $H^1(\Om)$.
Therefore, we can write that if\footnote{Under the condition~\eqref{C1tilde} and the fact that $3m-2\geq 1$ for all $m\geq 1$, in fact we always have $qN\geq 1$. In the same way we will note that $-qM\geq 1$ is satified in the $2$-dimensional case taking~$q$ large enough and in the $3$-dimensional case taking $q=6$ and using the assumptions given on page~\pageref{5} for~$\alpha$ and~$n$.} $qN\geq 1$ and $-qM\geq 1$ then
\begin{equation}\label{36}
  \vrh_\vep\rightarrow\vrh\ \ {\rm in}\ L^{qN}(\Om),
\end{equation}
\begin{equation}\label{37}
  \vrh_\vep^{-1}\rightarrow\vrh^{-1}\ \ {\rm in}\ L^{-qM}(\Om).
\end{equation}
Refering to the conditions \eqref{6} and \eqref{7}, we also conclude that, for all $q$ also satisfying\footnote{As previously, these two conditions are satisfied, using the condition~\eqref{C1tilde}, the fact that $3m-2 \geq m$ for all $m\geq 1$ and $\frac{-2qM}{n}\geq 2$ (since we previously prove that $-qM\geq 1$ and since $n<1$).} $2qN\geq m$ and $-2qM\geq n$
\begin{equation}\label{38}
  \sqrt{\mu(\vrh_\vep)}\rightarrow\sqrt{\mu(\vrh)}\ \ {\rm in}\ L^\frac{2qN}{m}(\Om),
\end{equation}
\begin{equation}\label{39}
  \frac{1}{\sqrt{\mu(\vrh_\vep)}}\rightarrow\frac{1}{\sqrt{\mu(\vrh)}}\ \ {\rm in}\ L^\frac{-2qM}{n}(\Om).
\end{equation}

\SUBSECTION{5.3}{Compactness on the velocity}

We know by \eqref{26} that $\sqrt{\mu(\vrh_\vep)}\pa_Z v_\vep$ is bounded in $L^2(\Om)$ and thus weakly converge in $L^2(\Om)$ to some~$f$.
>From the identity
\[
\pa_Z v_\vep=\frac{1}{\sqrt{\mu(\vrh_\vep)}}\sqrt{\mu(\vrh_\vep)}\pa_Z v_\vep,
\]
we also conclude that $\pa_Z v_\vep$ is bounded in $L^r(\Om)$ with $\frac{1}{r}=\frac{1}{2}-\frac{n}{2qM}$.\par
We note that $r \geq 1$ since we have previously proved that $-qM \geq 1 > n$.
We have
\begin{equation}\label{40}
\pa_Z v_\vep\rightharpoonup \pa_Z v\ \ {\rm in}\ L_w^r(\Om).
\end{equation}
Remark that we necessarily have $r<2$.\par
As for the derivatives with respect to $Z$, we know, using the bound \eqref{25} and those which come from the convergence \eqref{39}, that $\vep\pa_x v_\vep$ is bounded in $L^r(\Om)$ and thus weakly converges.
Thus, using the Poincar\'e inequality, we also get the bound of $v_\vep$ in $L^r(\Om)$ and then
\begin{equation}\label{41}
  v_\vep\rightharpoonup v\ {\rm in}\ L^r_w(\Om),
\end{equation}
\begin{equation}\label{42}
  \vep v_\vep\rightharpoonup 0\ {\rm in}\ W^{1,r}_w(\Om).
\end{equation}
Taking the same way, we also have
\begin{equation}\label{43}
  \vep w_\vep\rightharpoonup w\ {\rm in}\ L^r_w(\Om),
\end{equation}
\begin{equation}\label{44}
  \vep\pa_Z w_\vep\rightharpoonup \pa_Z w\ {\rm in}\ L^r_w(\Om),
\end{equation}
\begin{equation}\label{45}
  \vep^2 w_\vep\rightharpoonup 0\ {\rm in}\ W^{1,r}_w(\Om).
\end{equation}
Moreover, thanks to the compactness $W^{1,r}(\Om)\subset L^s(\Om)$ for all $s < r'$ where $\frac{1}{r'}=\frac{1}{r}-\frac{1}{d}$ we can write
\begin{equation}\label{46}
  \vep v_\vep\rightarrow 0\ {\rm in}\ L^s(\Om),\ \forall s < \frac{rd}{d-r},
\end{equation}
\begin{equation}\label{47}
  \vep^2 w_\vep\rightarrow 0\ {\rm in}\ L^s(\Om),\ \forall s < \frac{rd}{d-r}.
\end{equation}
%

\SUBSECTION{5.4}{Limit in the momentum equation}

Let us rewrite the two components~\eqref{18} and~\eqref{20} of the momentum equation:
\begin{eqnarray}
\vep^2\vrh_\vep v_\vep\pa_x v_\vep +\vep^2 \vrh_\vep w_\vep \pa_Z v_\vep & = & 2\vep^2\pa_x(\mu(\vrh_\vep) \pa_x v_\vep) + \boldsymbol{ \pa_Z(\mu(\vrh_\vep) \pa_Z v_\vep) } +\vep^2 \pa_Z(\mu(\vrh_\vep) \pa_x w_\vep) \nonumber \\
 & & + \vep^2\pa_x(\lambda(\vrh_\vep) (\pa_x v_\vep + \pa_Z w_\vep)) \boldsymbol{ - \pa_x(P(\vrh_\vep)) } + r_0\vep^2 \vrh_\vep (v_\vep^2+\vep^2 w_\vep^2)^{\frac{1}{2}} v_\vep, \nonumber \\
\vep^4\vrh_\vep v_\vep\pa_x w_\vep + \vep^4\vrh_\vep w_\vep \pa_Z w_\vep & = & \vep^2\pa_x(\mu(\vrh_\vep)\pa_Z v_\vep) + \vep^4\pa_x(\mu(\vrh_\vep) \pa_x w_\vep) + 2\vep^2\pa_Z(\mu(\vrh_\vep) \pa_Z w_\vep) \nonumber \\
 & & + \vep^2 \pa_Z(\lambda(\vrh_\vep) (\pa_x v_\vep + \pa_Z w_\vep)) \boldsymbol{ - \pa_Z(P(\vrh_\vep)) } + r_0\vep^4 \vrh_\vep (v_\vep^2+\vep^2 w_\vep^2)^{\frac{1}{2}}w_\vep. \nonumber
\end{eqnarray}
We initially will show that the bold terms admit limits when~$\eps$ tends to~$0$, then that all the other terms tend to zero.
The method is exactly the same one as that developed in part~\ssubsectionref{3.4.3}.
We just will specify the dependences in the parameter~$\eps$.
\par
$\bullet$ For instance, putting together assumptions on the pressure \eqref{9}--\eqref{10}, the estimates \eqref{29}--\eqref{32} and the strong convergences of the density \eqref{36}--\eqref{37}, one obtains the convergence of~$\dx(P(\vrh_\vep))$ to~$\dx(P(\vrh))$ and the convergence of~$\dZ(P(\vrh_\vep))$ to~$\dZ(P(\vrh))$ in the sense of distributions on~$\Om$, since we introduce the same hypotheses on the coefficients, see for instance the convergence of the term~$T_1$ in subsection~\ssubsectionref{3.4.3} and the assumption~\eqref{C2}.

\par
$\bullet$ As for the term $T_3$ in the subsection~\ssubsectionref{3.4.3}, under the condition~\eqref{C'3} we have
\begin{equation}\label{48}
  \sqrt{\mu(\vrh_\vep)}\pa_Z v_\vep\rightharpoonup \sqrt{\mu(\vrh)}\pa_Z v\ \ {\rm in}\ L^2(\Om).
\end{equation}
Since $m<qN$, we end to the convergence:
\begin{equation}\label{49}
  \mu(\vrh_\vep)\pa_Z v_\vep\rightarrow \mu(\vrh)\pa_Z v\ \ {\rm in}\ \cD'(\Om).
\end{equation}
We are now going to show that the other terms tend to~$0$ in the sense of distributions when $\vep$ goes to~$0$.\par
%
%
%
%
%
%
$\bullet$ As for \eqref{48} and \eqref{49}, we also get, for the $x$-derivatives,
\begin{equation}\label{50}
\vep \sqrt{\mu(\vrh_\vep)}\pa_x v_\vep\rightharpoonup 0\ \ {\rm in}\ L^2(\Om),
\end{equation}
\begin{equation}\label{51}
  \vep \mu(\vrh_\vep)\pa_x v_\vep\rightarrow 0\ \ {\rm in}\ \cD'(\Om).
\end{equation}
Therefore, we obtain the convergence of $2\vep^2\pa_x(\mu(\vrh_\vep) \pa_x v_\vep)\rightarrow 0$ in $\cD'(\Om)$ and also $\vep^2\pa_x(\lambda(\vrh_\vep)\pa_x v_\vep) \rightarrow 0$ and $\vep^2 \pa_Z(\lambda(\vrh_\vep)\pa_x v_\vep)\rightarrow 0$ since $\la(s)=2(s\mu'(s)-\mu(s))$ satisfies the same integrabilities as $\mu$.
Taking the same strategy as for the convergences already obtained through \eqref{48}, \eqref{49}, \eqref{50} and \eqref{51}, we finally can give the ones of the vertical velocity $w_\vep$ :
\begin{eqnarray}
     \vep \sqrt{\mu(\vrh_\vep)}\pa_Z w_\vep \rightharpoonup & \sqrt{\mu(\vrh)}\pa_Z w\ & \ {\rm in}\ L^2(\Om), \label{52} \\
     \vep \mu(\vrh_\vep)\pa_Z w_\vep \rightarrow & \mu(\vrh)\pa_Z w\ & \ {\rm in}\ \cD'(\Om), \label{53} \\
     \vep^2 \sqrt{\mu(\vrh_\vep)}\pa_x w_\vep \rightharpoonup & 0\ & \ {\rm in}\ L^2(\Om), \label{54} \\
     \vep^2 \mu(\vrh_\vep)\pa_x w_\vep \rightarrow & 0\ & \ {\rm in}\ \cD'(\Om), \label{55}
\end{eqnarray}
which answer the questions of convergence for the viscous terms containing $w_\vep$.\par
$\bullet$ For the friction terms, we will just do it for the term $r_0\eps^2 \rho_\eps v_\eps^2$ (the other terms may be treated using exactly the same way, moreover the convective terms are also treated in the same way).
Using the same method as in Subsection~\ssubsectionref{3.4.3}, in particular when we have treated the term~$T_2$, we have
\begin{equation*}
  r_0\eps^2 \rho_\eps v_\eps^2  = r_0 \eps^{\frac{2}{3}} \, (\eps² \rho_\eps v_\eps^3)^{\frac{2}{3}} \, \rho_\eps^{\frac{1}{3}}.
\end{equation*}
>From the strong convergence of~$\rho_\eps$ and the weak convergence of $\eps² \rho_\eps v_\eps^3$ we deduce that $r_0\eps^2 \rho_\eps v_\eps^2$ tends to~$0$.
\noindent
Thus, we can insure that \eqref{18}--\eqref{19} converge to \eqref{21}--\eqref{22} in the sense of distributions.

\SUBSECTION{5.5}{Limit in the mass equation}

To pass to the limit $\eps\rightarrow 0$ in the mass equation
$
\pa_x (\vrh_\vep v_\vep) + \pa_Z (\vrh_\vep w_\vep) = 0,
$
the main difficulty comes from the fact that the vertical velocity~$w_\eps$ does not have a limit. We thus will use the following equivalent form
\begin{equation*}
\pa_x \Big( \int_0^{h} \vrh_\vep v_\vep\, dZ \Big) = 0.
\end{equation*}

It is an equivalent form in the following sens: if a velocity~$v_\eps$ such that $\dsp \pa_x \Big( \int_0^{h} \vrh_\vep v_\vep\, dZ \Big) = 0$ exists then we can build a vertical velocity $w_\eps$ such that $\pa_x (\vrh_\vep v_\vep) + \pa_Z (\vrh_\vep w_\vep) = 0$ and $w_\eps|_{Z=0}=w_\eps|_{Z=h}=0$.
This is enough to define $\dsp w_\eps = - \frac{1}{\rho} \int_0^Z \pa_x (\vrh_\vep v_\vep)$.
\par
To prove that $\rho_\eps v_\eps$ tends to $\rho v$ we write
\begin{equation*}
  \rho_\eps v_\eps = (\rho_\eps v_\eps^3)^{\frac{1}{3}} \, \rho_\eps^{\frac{2}{3}}.
\end{equation*}
Using the strong convergence of~$\rho_\eps$ and the bound on $\rho_k v_k^3$ in $L^1(\Omega)$, we obtain (as soon as $qN\geq 1$, which is implied by the condition~\eqref{C1})
\begin{equation*}
  \rho_\eps v_\eps \rightharpoonup g\, \rho^{\frac{2}{3}} \ \ {\rm in}\ L^{1}(\Om),
\end{equation*}
where $g$ is the weak limit of $\big( \rho_\vep |\u_\vep|^3 \big)^{\frac{1}{3}}$ in $L^{3}(\Omega)$.
To identify the limit~$g$, we  use the strong convergence for the density and the weak convergence for the velocity:
\begin{equation*}
  \rho_\vep^{\frac{1}{3}} \to \rho^{\frac{1}{3}} \ \ {\rm in}\ L^{3qN}(\Om)
\qquad \text{and} \qquad
  v_\vep \rightharpoonup v \ \ {\rm in}\ L^r(\Om).
\end{equation*}
We deduce that $g=\rho^{\frac{1}{3}}v$ if $\frac{1}{r}+\frac{1}{3qN}\leq 1$, condition which is implied by the condition~\eqref{C'3}.
\par
We deduce the following convergence
\begin{equation*}
\pa_x \Big( \int_0^{h} \vrh_\vep v_\vep\, dZ \Big) 
\rightarrow
\pa_x \Big( \int_0^{h} \vrh v\, dZ \Big) \quad \text{in $\cD'(\Om)$.}
\end{equation*}

\SUBSECTION{5.6}{Pressure and viscosity conditions}

In addition to the conditions given in Subsection \subsectionref{3.4}, we have to suppose that (see assumption~\eqref{C3})
\begin{equation*}
  qM \leq 4-3m.
\end{equation*}
{\bf In the two dimensional case}, $q$ can be chosen as large as we need, thus many inequalities are satisfied and Theorem \ref{th-convergence} holds only with conditions \eqref{conditions-coeff0}.
\par
{\bf In the three dimensional case}, $q$ is any number smaller than $6$.
The additive conditions on the pressure and viscosities coefficients can be summarized by:
\begin{equation*}
m < \al-n+\frac{7}{3}.
\end{equation*}
This condition exactly corresponds to the condition~\eqref{condition+}.
\par
\vskip 20pt
\hrule
\vskip 10pt

{\it The authors would like to thank Didier Bresch for his useful devices about that kind of studies.
The first author of this work has been partially supported by the ANR project n° ANR-08-JCJC-0104-01~: RUGO (Analyse et calcul des effets de rugosités sur les écoulements).}

\vskip 10pt
\hrule
\vskip 20pt

\BIB

\bib{Arregui}
\sc I. Arregui, J. Cendán, C. Parés, C. Vázquez.
\rm Optimization of a duality method for the compressible Reynolds equation.
\it Numerical mathematics and advanced applications, 
\rm 319--327, Springer, Berlin, 2006.

\bib{Bayada-Chambat}
\sc G. Bayada, M. Chambat.
\rm The transition between the Stokes equations and the Reynolds equation: a mathematical proof.
\it Appl. Math. Optim. 
\rm 14  (1986),  no. 1, 73--93.

\bib{Bhushan}
\sc B. Bhushan, K. Tonder.
\rm Roughness-induced shear-and squeeze-film effects in magnetic recording-Part I: Analysis.
\it ASME J. Tribol. 
\rm 111 (1989), pp. 220--227.

\bib{BD1}
\sc D. Bresch, B. Desjardins, 
\rm Existence of global weak solutions for a 2D viscous shallow water equations and convergence to the quasigeostrophic
model, 
\it Comm. Math. Phys. 
\rm 238 (1-2) (2003), 211-223.

\bib{BD2}
\sc D. Bresch, B. Desjardins, 
\rm Some diffusive capillary models of Korteweg type, C. R. Acad. Sci., Paris, Section mécanique 332 (11) (2004), 881-886.

\bib{BD3}
\sc D. Bresch, B. Desjardins, 
\rm On the construction of approximate solutions for 2D viscous shallow water model and for compressible Navier-Stokes models, J. Math. Pures Appl. 86 (4) (2006) 362-368.

\bib{BD}
\sc D. Bresch, B. Desjardins,
\rm On the existence of global weak solutions to the Navier-Stokes equations for viscous compressible and heat conducting fluids.
\it J. Math. Pures et Appliqu\'ees,
\rm (9) 87, (2007), no.1, 57--90.

\bib{BDG}
\sc D. Bresch, B. Desjardins, D. Gerard-Varet,
\rm On compressible Navier-Stokes equations with density dependent viscosities in bounded domains.
\it 
\rm 

\bib{BDL}
\sc D. Bresch, B. Desjardins, C.K. Lin, 
\rm On some compressible fluid models: Korteweg, lubrication and shallow water systems.
\it Comm. Partial Diff. Eqs,
\rm {28}, No. 3-4, (2003), 1009--1037.

\bib{Jai}
\sc G. Buscaglia, I. Ciuperca, M. Jai,
\rm Existence and uniqueness for several non-linear elliptic problems arising in lubrication theory.
\it J. Differential Equations
\rm 218 (2005), no. 1, 187--215.

\bib{Buscaglia}
\sc G. Buscaglia, M. Jai,
\rm A new numerical scheme for non uniform homogenized problems: application to the non linear Reynolds compressible equation.
\it Math. Probl. Eng. 
\rm 7 (2001), no. 4, 355--378.

\bib{Crone}
\sc R.M. Crone, M.S. Jhon, B. Bhushan, T.E. Karis.
\rm Modeling the flying characteristics of a rough magnetic head over a rough rigid-disk surface.
\it ASME J. Tribol. 
\rm 113 (1991), pp. 739--749.

\bib{FNP}
\sc E. Feireisl, A. Novotny, H. Petzeltova, 
\rm On the existence of globally defined weak solutions to the Navier-Stokes equations of compressible
isentropic fluids, J. Math. Fluid Dynam. 3 (2001) 358-392.

\bib{F}
\sc E. Feireisl, 
\it Dynamics of Viscous Compressible Fluids, 
\rm Oxford Science Publication, Oxford, 2004.

\bib{L}
\sc P.-L. Lions, 
\rm Compacité des solutions des équations de Navier-Stokes compressibles isentropiques, 
\it C. R. Acad. Sci., Paris, Sér. 
\rm I 317 (1993) 115-120.

\bib{L2}
\sc P.-L. Lions, 
\it Mathematical Topics in Fluid Dynamics, vol. 2, 
\rm Compressible Models, Oxford Science Publication, Oxford, 1998.

\bib{M}
\sc S. Martin, 
\rm Contribution à la modélisation de phénomènes de frontière libre en mécanique des films minces, 
\rm Thèse de doctorat de l'Institut National des Sciences Appliquées de de Lyon (2005).

\bib{MarPal-Sta}
\sc E. Marusic-Paloka, M. Starcevic,
\rm Rigorous justification of the Reynolds equations for gas lubrication.
\it Comptes Rendus Mecanique, Paris,
\rm 333 (2005), 534-541. 

\bib{MarPal-Sta09}
\sc E. Marusic-Paloka, M. Starcevic,
\rm Derivation of Reynolds equation for gas lubrication via asymptotic analysis of the compressible Navier-Stokes system.
\it Nonlinear Analysis: Real World Applications, 
\rm (2009).

\bib{Mik}
\sc A. Mikelic,
\rm On the justification of the Reynolds equation, describing isentropic compressible flows trough a tiny pore.
\it Ann. Univ. Ferrara Sez. VII Sci. Mat. 
\rm 53 (2007), no. 1, 95--106. 

\bib{Nazarov}
\sc S.A. Nazarov, Y.G. Videman.
\rm An improved nonlinear Reynolds equation for a thin flow of a viscous incompressible fluid.
\it Vestnik St. Petersburg Univ. Math.  
\rm 41 (2008),  no. 2, 171--175.

\bib{NS}
\sc A. Novotny, I. Straskraba, 
\it Introduction to the Mathematical Theory of Compressible Flow, 
\rm Oxford Science Publication, Oxford, 2004.

\bib{Reynolds}
\sc O. Reynolds.
\rm On the theory of lubrication and its application to M.~Beauchamp Tower's experiments.
\it Phil. Trans. Roy. Soc. London,
\rm  A117, 157--234, 1886.

\bib{Sart}
\sc R. Sart.
\rm Viscous Augmented Born-Infeld for magnetohydrodynamic flows.
\it J. Math. Fluid Mech.,
\rm  11, (2008), 1--25.

\ENDBIB

\bye


\end{document}